\documentclass[10pt]{article}

\usepackage{jmlr2e}

\usepackage{amsmath}
\usepackage{bm}
\usepackage{amsfonts}
\usepackage{algorithm}
\usepackage{algorithmic}

\usepackage{wrapfig,epsfig}

\usepackage{graphicx}

\usepackage{enumitem}

\usepackage{xcolor}

\usepackage{authoraftertitle}

\usepackage{grffile}

\newtheorem{Definition}{Definition}
\newtheorem{Lemma}{Lemma}
\newtheorem{Theorem}{Theorem}
\newtheorem{Assumption}{Assumption}
\newtheorem{Remark}{Remark}
\newtheorem{Corollary}{Corollary}

\newcommand{\AlgorithmFullName}{Byzantine-resilient decentralized TD($\lambda$) with linear function approximation}

\title{Byzantine-Resilient Decentralized TD Learning with\\ Linear Function Approximation}

\author{
    Zhaoxian Wu$^{\star}$ \qquad 
    Han Shen$^{\dag}$ \qquad  
    Tianyi Chen$^{\dag}$ \qquad 
    Qing Ling$^{\star}$\\\\
    {\centering
        \addr 
        $^{\star}$\textit{Sun Yat-Sen University, Guangzhou, Guangdong 510006, China} \\
        $^{\dag}$\textit{Rensselaer Polytechnic Institute, Troy, New York 12180, USA}\\
        \email \texttt{wuzhx23@mail2.sysu.edu.cn}~~
            \{\texttt{shenh5,chent18\}@rpi.edu}~~
            \texttt{lingqing556@mail.sysu.edu.cn}
    }
}

\ShortHeadings{Byzantine-Resilient Decentralized TD Learning with\ Linear Function Approximation}{Wu, Shen, Chen and Ling}
\firstpageno{1}

\begin{document}

\maketitle
\thispagestyle{empty}
\begin{abstract}
	This paper considers the policy evaluation problem in a multi-agent reinforcement learning (MARL) environment over decentralized and directed networks. The focus is on decentralized temporal difference (TD) learning with linear function approximation in the presence of unreliable or even malicious agents, termed as Byzantine agents. In order to evaluate the quality of a fixed policy in a common environment, agents usually run decentralized TD($\lambda$) collaboratively. However, when some Byzantine agents behave adversarially, decentralized TD($\lambda$) is unable to learn an accurate linear approximation for the true value function. We propose a trimmed-mean based Byzantine-resilient decentralized TD($\lambda$) algorithm to perform policy evaluation in this setting. We establish the finite-time convergence rate, as well as the asymptotic learning error in the presence of Byzantine agents. Numerical experiments corroborate the robustness of the proposed algorithm.
	
	
	
\end{abstract}

\section{Introduction}

Reinforcement learning (RL) is a promising paradigm of modern artificial intelligence. An agent in RL acts under a specific policy, interacts with an unknown environment, observes environment state changes, receives rewards, and improves the policy. RL has demonstrated great potential in various applications, including autonomous driving \cite{chen2015deepdriving}, robotics \cite{gu2017deep}, power networks \cite{kar2012qd}, to name a few.



Among many algorithms that are central to RL is the temporal-difference (TD) learning approach.
TD learning provides a unified framework to evaluate a policy in terms of the long-term discounted accumulative reward. Within the family of TD algorithms, TD($\lambda$) is a popular one \cite{sutton1988learning}. However, the classical TD($\lambda$) is a tabular-based approach, which stores the entry-wise value function state-by-state. In many RL applications, the number of elements in the state space is large. Thus, it is impractical or even impossible to evaluate the value function on a per state basis.
This is also referred to as the curse of dimensionality. Therefore, TD with function approximation is usually preferable, which includes linear or non-linear approximators \cite{mnih2015human, silver2016mastering}.


Going beyond the single-agent RL setting, many practical applications contain rewards distributed over multiple agents without a central coordinator. In this multi-agent RL (MARL) scenario, agents cooperate with each other to explore an unknown environment and accomplish a specific task  \cite{zhang2018fully,wai2018multi,chen2018communication,lee2020optimization,Zhang2019multi}. A broad range of applications over mobile sensor networks \cite{cortes2004coverage} and power networks \cite{kar2012qd} can be modeled as MARL problems. Typically, agents in MARL are assumed to be reliable; that is, they share accurate information with others. In this circumstance, TD($\lambda$) can be implemented over a MARL environment in a decentralized manner \cite{doan2019finite}.


However, the agents in MARL are not always reliable. Communication errors, data corruptions or even malicious attacks may happen during the training process. We model all these unreliable behaviors as \textit{Byzantine attacks} \cite{lamport1982,Bajwa2019-review}. In this Byzantine attack model, adversarial agents are able to arbitrarily manipulate their outputs, collude with each other and inject false information into the MARL environment. Decentralized TD, and policy evaluation in general, in the presence of Byzantine agents are yet less understood areas.


To perform decentralized policy evaluation under Byzantine attacks, we propose a novel trimmed mean-based \AlgorithmFullName{} algorithm that we term Byrd-TD($\lambda$). Different from existing decentralized TD($\lambda$) such as that in \cite{doan2019finite}, Byrd-TD($\lambda$) is a robust version of TD($\lambda$). It uses a robust aggregation technique to alleviate the negative effect of Byzantine agents, and guarantees global convergence even under Byzantine attacks. However, the theoretical analysis, especially the non-asymptotic analysis of Byrd-TD($\lambda$), is nontrivial. First, the existence of Byzantine agents prevents the honest agents from obtaining accurate information from their neighbors, and the induced error propagates through the TD($\lambda$) dynamics. Second, the fact that the TD($\lambda$) update does not follow stochastic gradient of any objective function brings a series of difficulties to the non-asymptotic analysis of Byrd-TD($\lambda$) by leveraging the first-order optimization toolbox.



\subsection{Related works}
TD learning plays an important role in policy evaluation for RL. The family of TD methods include gradient TD \cite{sutton2008convergent,sutton2009fast}, least-squares TD \cite{lazaric2012finite,bradtke1996linear}, least-squares policy evaluation \cite{nedic2003least}. As a representative algorithm, TD($\lambda$)  provides a unified framework for policy evaluation \cite{sutton1988learning}, with the parameter $\lambda \in [0,1]$ controlling the trade-off between approximation accuracy and convergence rate. Several existing works have shown the convergence of TD($\lambda$) with linear function approximation.
The works of \cite{tsitsiklis1997analysis, dayan1992convergence} focus on the asymptotic convergence and other steady-state properties of TD($\lambda$). Several works provide finite-time analysis \cite{bhandari2018finite, lakshminarayanan2018linear, dalal2018finite,sun2020adaptive}, but their results hold only when using projection operation or assuming independent and identically distributed (i.i.d.) noise. Some of these assumptions have been relaxed in the recent works of \cite{srikant2019finite}. Building upon these works, convergence analysis of decentralized TD has been studied recently in \cite{doan2019finite,sun2020finite,cassano2020multi}.

Some exploratory efforts have been made to robustify RL approaches. In the single-agent setting, algorithms have been developed to alleviate the negative effect from approximation error \cite{lu2019family}, or to robustify the learned policy in making decisions \cite{tirinzoni2018policyconditioned,Petrik2019beyond}. However, little is known on how to robustify RL methods during the training stage in the single-agent RL setting, and in the MARL setting -- an even less explored area.

On the other spectrum, designing Byzantine-resilient algorithms has become a popular topic in distributed and decentralized machine learning, with applications mainly in supervised learning. Most algorithms in this regime replace the mean or weighted mean aggregation rules in Byzantine-free stochastic gradient methods with robust aggregation rules, such as geometric median \cite{chen2019DistributedSM, xie2018GeneralizedBS}, coordinate-wise median \cite{yin2018ByzantineRobustDL}, trimmed mean \cite{xie2018phocas}, Krum \cite{blanchard2017MachineLW}, Bulyan \cite{guerraoui2018hidden}, and RSA \cite{li2019RSABS}, or robustify the updates with malicious agent identification \cite{alistarh2018Byzantine}. These aforementioned methods require a central node to coordinate the agents, which is different to our decentralized MARL setting. Existing decentralized Byzantine-resilient algorithms include the subgradient method with total variation-norm penalization \cite{ben2015robust, peng2020byzantine}, BRIDGE \cite{yang2019bridge}, and MOZI \cite{guo2020towards}. However, these algorithms are all targeted to supervised learning either in i.i.d. or deterministic settings.
Policy evaluation in MARL with adversarial agents is, by its nature, a stochastic optimization problem with non-i.i.d. and decentralized data. To the best of our knowledge, no Byzantine-resilient SGD-based algorithm has been developed in this regime.
Therefore, existing Byzantine-resilient algorithms can not be applicable to policy evaluation, especially in MARL.

\begin{figure}[t]
	\centerline{\includegraphics[width=0.6\columnwidth]{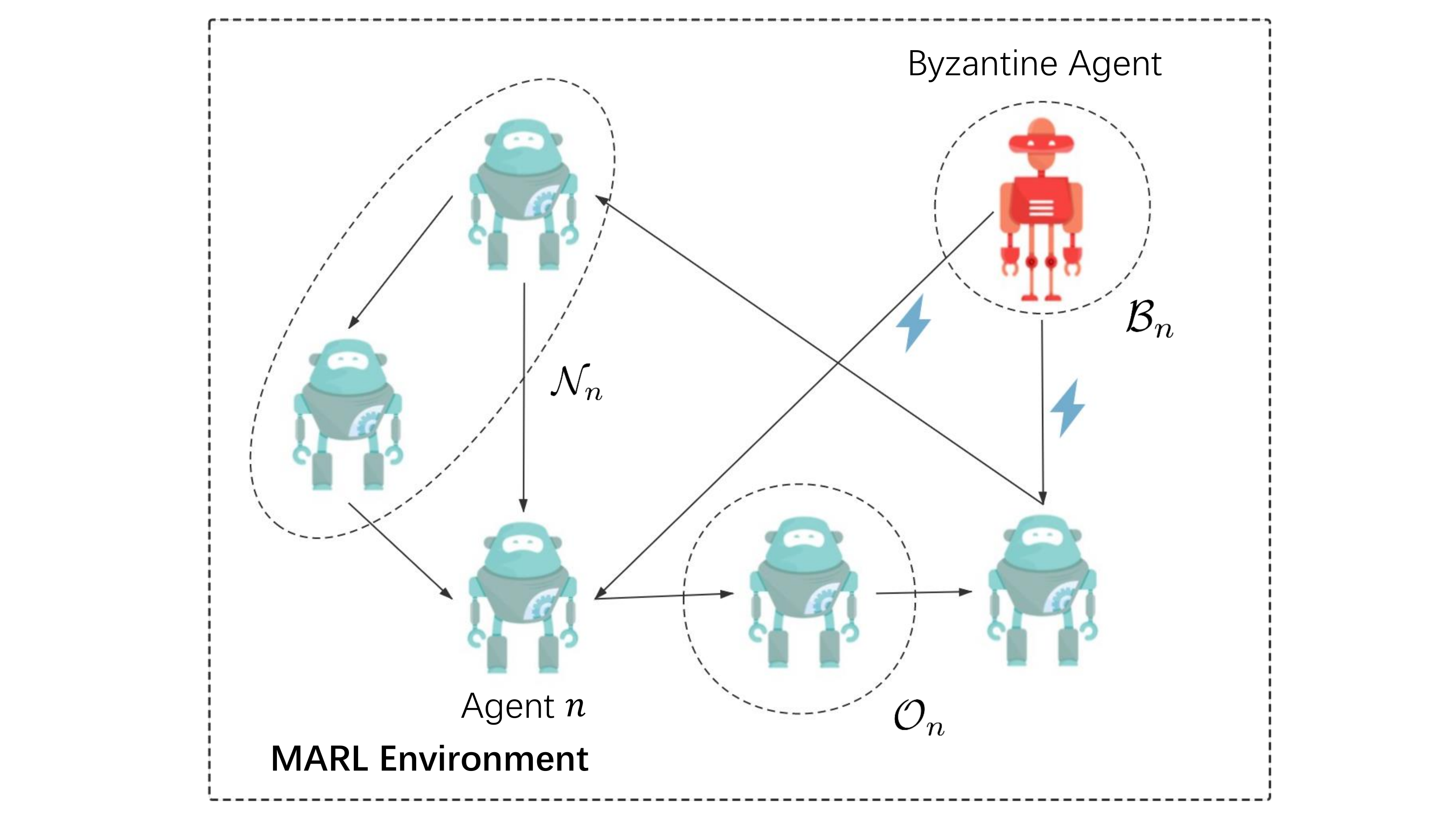}}
	\caption{MARL with Byzantine agents.}
	\label{fig:ByrdTD-diag}
\end{figure}


In the Byzantine-resilient decentralized TD($\lambda$) algorithm, each agent uses trimmed mean to aggregate neighboring messages. This strategy naturally generates a directed and time-varying graph that is a subset of the underlying communication network. Existing works have investigated decentralized optimization when the underlying communication network is directed and time-varying \cite{yuan2014randomized,nedic2017achieving,nedic2014distributed}. However, to the best of our knowledge, our work is the first to analyze decentralized TD($\lambda$) over a directed and time-varying graph.

\subsection{Our contributions} The contributions of this paper are summarized as follows.

1. Targeting the robust policy evaluation in MARL, we propose a \AlgorithmFullName{}, abbreviated as Byrd-TD($\lambda$), over decentralized and directed networks. Byrd-TD($\lambda$) generalizes the robust aggregation rule \cite{yang2019bridge,vaidya2012matrix} to the regime of stochastic optimization, which is of independent interest.

2. To analyze our algorithm, we generalize the existing trimmed mean analysis to the stochastic setting with time-varying graphs. We establish the non-asymptotic analysis for Byrd-TD($\lambda$) under Markovian observations.
Our analysis implies that the network topology influences the convergence rate, and the asymptotic  error is determined by the variation of local rewards and the degree of network unsaturation.


3. We conduct experiments using a standard MARL benchmark, the cooperative navigation task, over various networks. It is confirmed that Byrd-TD($\lambda$) can successfully learn the linear value function estimate under Byzantine attacks.


%

\section{Byzantine-resilient Decentralized TD}
In this section, we review some background on the policy evaluation, and then introduce our algorithm.

\subsection{Preliminaries on policy evaluation and TD}
Consider an MARL environment consisting of $N=|\mathcal{N}|$ honest agents and $B=|\mathcal{B}|$ Byzantine agents, where $\mathcal{N}$ and $\mathcal{B}$ denote the sets of honest agents and Byzantine agents, respectively. Suppose that the agents can only communicate with each other through a given directed network $\mathcal{G} := (\mathcal{N}\cup\mathcal{B}, \mathcal{E})$, where $\mathcal{E} \subseteq (\mathcal{N}\cup\mathcal{B})\times(\mathcal{N}\cup\mathcal{B})$ is the edge set without self-link. An honest agent $n$ can only send information to its out-neighbors in $\mathcal{O}_n = \{m | (n, m)\in\mathcal{E}\}$ and receive information from its in-neighbors in $\mathcal{N}_n\cup\mathcal{B}_n$, where $\mathcal{N}_n=\{m | (m, n)\in\mathcal{E}, m\in\mathcal{N}\}$ and $\mathcal{B}_n=\{m | (m, n)\in\mathcal{E}, m\in\mathcal{B}\}$ are the honest (excluding $n$ itself) and Byzantine in-neighbor sets of honest agent $n$, respectively. We define the numbers of its honest and Byzantine in-neighbors by $N_n=|\mathcal{N}_n|$ and $B_n=|\mathcal{B}_n|$, respectively. We illustrate MARL with Byzantine agents in Figure 1. In this environment, our goal is to evaluate a fixed policy in a Markov reward process (MRP).



A typical MRP \footnote{To simplify the notation, we will drop the dependence on the policy as those considered in the Markov decision process, since the policy is fixed.} can be described by a 5-tuple  $(\mathcal{S}, \mathcal{P}, \gamma, \rho^0, \{\mathcal{R}_n\}_{n\in{\cal N}})$, where $\mathcal{S}$ is a set of states shared by all agents, $\mathcal{P}$ is the space of the state transition kernels with the probability of transiting from state $s$ to $s'$ being defined as $\mathcal{P}(s, s')$, $\gamma\in(0, 1)$ is the discounting factor, $\rho^0$ is the initial state distribution, and $\mathcal{R}_n(s, s')$ is the local reward of agent $n$ from state $s$ to $s'$. At each state transition $s\rightarrow s'$, every honest agent $n$ only observes its local reward $\mathcal{R}_n(s, s')$ from the environment and has access to its local reward function $\mathcal{R}_n(s, s')$. The global reward function is defined as the average of the local rewards $\frac{1}{N}\sum_{n\in\mathcal{N}} \mathcal{R}_n(s, s')$. The goal of policy evaluation is to find the discounted accumulative reward starting from the initial state $s$, which is referred to as the value function of state $s$, given by
\begin{equation}
	\label{definition:valueFunction}
	\mathcal{V}(s) := \mathbb{E}
	\left[
	\left.
	\frac{1}{N}\sum_{k=0}^{\infty}\gamma^k \sum_{n\in\mathcal{N}} \mathcal{R}_n(s^k, s^{k+1})
	\right| s^0 = s
	\right],
\end{equation}
where the expectation is taken over all possible state trajectories $\{s^0, s^1, \dots, s^k, \dots\}$.

\textbf{TD($\lambda$) with function approximation.} One way to find the value function is to solve the so-termed Bellman equation, which is a system of linear equations of $\mathcal{V}(s)$, given by
\begin{equation}
	\label{equation:BellmanEqu}
	\mathcal{V}(s)
	= \frac{1}{N}\sum_{s'\in\mathcal{S}}
	\mathcal{P}(s, s')
	\Big(
	\sum_{n\in\mathcal{N}} \mathcal{R}_n(s, s') + \gamma\mathcal{V}(s')
	\Big),~\forall s\in{\cal S}.
\end{equation}
If $\mathcal{P}(s, s')$ is known and the number of elements in the state space $|\mathcal{S}|$ is small, one can find $\mathcal{V}(s)$ by solving this system of linear equations. On the other hand, if $\mathcal{P}(s, s')$ is unknown yet $|\mathcal{S}|$ is manageable, tabular-based TD methods have impressive performance.
Unfortunately, when the number of elements in $\mathcal{S}$ is large or even infinite, it is impractical to find an exact solution or even evaluate a solution to \eqref{equation:BellmanEqu}. Instead, one has to resort to function approximation of $\mathcal{V}(s)$.

Consider an approximation function $V_{\bm{\theta}}(s)$ parameterized by a vector $\bm{\theta}\in\mathbb{R}^D$.
In this paper, we consider linear approximation, such that $V_{\bm{\theta}}(s) := \bm{\phi}(s)^{\top} \bm{\theta}$, where $\bm{\phi}(s)\in\mathbb{R}^D$ is the pre-defined feature vector of state $s$. For the sake of computational tractability, we often have the number of unknown parameters $D\ll |\cal S|$. With the above linear approximation, the task is to search for an appropriate parameter $\bm{\theta}$ that solves the following optimization problem
\begin{equation}
	\label{equation:problem}
	\min_{\bm{\theta}} ~ F(\bm{\theta}) := \frac{1}{2} \sum_{s\in\mathcal{S}} \rho(s)(V_{\bm{\theta}}(s) - \mathcal{V}(s))^2,
\end{equation}
where $\rho(s)$ is the stationary distribution of state $s$ associated with the MRP. Observe that $F(\bm{\theta})$ represents the weighted least-squares \textit{approximation error} with the parameter $\bm{\theta}$.


When there is a central node to coordinate the learning process, one can solve \eqref{equation:problem} using the  gradient descent update 
\begin{equation}
	{\bm{\theta}}^{k+1}={\bm{\theta}}^{k}-\eta^k \nabla F({\bm{\theta}}^k),
\end{equation}
where $\eta^k > 0$ is the step size at step $k$. However, the exact gradient $\nabla F({\bm{\theta}}^k)$ or even its unbiased stochastic estimate is unavailable in practice. This motivates the use of TD ($\lambda$).
Let us define the \textit{eligibility trace} as
\begin{equation}
	\label{definition:eligibilityTrace}
	\bm{z}^k
	:= \sum_{\kappa=0}^k(\gamma\lambda)^{k-\kappa}\nabla_{\bm{\theta}^\kappa} V_{\bm{\theta}^\kappa}(s)
	= \sum_{\kappa=0}^k(\gamma\lambda)^{k-\kappa}\bm{\phi}(s^\kappa),
\end{equation}
where $\lambda \in [0,1]$ is a predefined parameter. Then, the update of TD($\lambda$) with linear function approximation can be written as
\begin{equation}
	\label{equation:update-TD-normal}
	{\bm{\theta}}^{k+1}={\bm{\theta}}^{k}
	+\frac{\eta^k}{N}\sum_{n\in\mathcal{N}} \left(
	r_n^k+\left(\gamma\bm{\phi}(s^{k+1})-\bm{\phi}(s^{k})\right)^{\top}{\bm{\theta}}^{k}
	\right) \bm{z}^k,
\end{equation}
where $r_n^k:=\mathcal{R}_n(s^k, s^{k+1})$. It is known that TD($\lambda$) reduces to TD($0$) when $\lambda=0$, and reduces to the Monte Carlo sampling method when $\lambda=1$ \cite[Chap. 7]{sutton2018reinforcement}.




\subsection{Decentralized TD($\lambda$) under Byzantine attacks}
\label{subsection:definition-trimmedmean}
To implement TD($\lambda$) in the MARL setting, every honest agent $n$ maintains its local copy $\bm{\theta}_n^k \in\mathbb{R}^D$. At step $k$, the honest agent $n$ observes the transition from state $s^k$ to $s^{k+1}$, calculates the local reward $r_n^k$, updates the eligibility trace by
\begin{equation}
	\label{eq:zzzzz}
	\bm{z}^k=\gamma\lambda\bm{z}^{k-1}+\bm{\phi}(s^k),
\end{equation}
and computes the local TD($\lambda$) increment at agent $n$ by
\begin{equation}
	\label{definition:TD_increment}
	\bm{g}_n^k(\bm{\theta}_n^k) := \left(r_n^k + (\gamma\bm{\phi}(s^{k+1})-\bm{\phi}(s^k))^{\top} \bm{\theta}_n^k \right)\bm{z}^k.
\end{equation}
After computing their local increments, the agents communicate with each other to exchange the latest local parameters. The honest agent $n$ collects $\bm{\theta}_m^k$ from all in-neighbors $m$ in $\mathcal{N}_n\cup\mathcal{B}_n$, averages them, and uses the local increment $\eta^k\bm{g}_n^k(\bm{\theta}_n^k)$ to update the local parameter as
\begin{equation}
	\label{equation:aggregate-decentralized-mean}
	\bm{\theta}_n^{k+1} = \frac{1}{N_n+B_n+1} \sum_{m\in\mathcal{N}_n\cup\mathcal{B}_n\cup\{n\}}\bm{\theta}_m^k
	+\eta^k\bm{g}_n^k(\bm{\theta}_n^k).
\end{equation}

\begin{wrapfigure}{r}{.55\linewidth}
\vspace*{-0.9cm}
\begin{minipage}{1\linewidth}
\begin{algorithm}[H]
	\caption{Byrd-TD($\lambda$)}
	\label{algorithm:TD_TrimmedMean}
	\begin{algorithmic}[0]
		\REQUIRE fixed policy; step size $\eta^k$; initial model $\bm{\theta}^0_n=\bm{\theta}^0$ and $\bm{z}^{-1}=\bm{0}$; $q_n$ for all honest agents $n$
		\FORALL {$k = 0, 1, 2, \cdots $}

		\FORALL {honest agents $n$}
		\STATE Send $\bm{\theta}_n^{k}$ to all out-neighbors $m\in\mathcal{O}_n$
		\STATE Receive $\bm{\theta}^{k}_m$ from all in-neighbors $m\in\mathcal{N}_n\cup\mathcal{B}_n$
		\ENDFOR

		\STATE Transition from state $s^k$ to $s^{k+1}$
		\FORALL {honest agents $n$}
		\STATE Observe the local reward $r_n^k$
		\STATE Compute $\bm{z}^k$ following \eqref{eq:zzzzz}
		\STATE Compute $\bm{g}_n^k(\bm{\theta}_n^k)$ following \eqref{definition:TD_increment}
		\FOR {$d=1, \dots, D$}
		\STATE Update $[\bm{\theta}_n^{k+1}]_d$ following \eqref{equation:aggregate-decentralized-trimmedMean}
		\ENDFOR
		\ENDFOR
		\ENDFOR
	\end{algorithmic}
\end{algorithm}
\end{minipage}
\vspace*{-0.3cm}
\end{wrapfigure}

The update rule \eqref{equation:aggregate-decentralized-mean} seems plausible. However, the main challenge of implementing \eqref{equation:aggregate-decentralized-mean} is that while $\bm{\theta}_m^{k}$ are true parameters for $m\in\mathcal{N}_n$, they can be arbitrary vectors for $m\in\mathcal{B}_n$. But the honest agent $n$ cannot distinguish Byzantine agents from all the in-neighbors. Motivated by the Byzantine-resilient supervised learning algorithms \cite{blanchard2017MachineLW,yang2019bridge,guo2020towards}, we will robustify \eqref{equation:aggregate-decentralized-mean} by incorporating \textit{coordinate-wise trimmed mean} to resist Byzantine attacks. Before taking averages over the messages received from their in-neighbors, the agents first discard the outliers on a per coordinate basis. To be specific, the honest agent $n$ estimates the number of Byzantine neighbors $B_n$ and chooses a trimming number $q_n \ge B_n$. At every dimension $d$, the honest agent $n$ discards the largest $q_n$ and smallest $q_n$ elements of the received messages before taking average. We use $\mathcal{N}_n^{k+}(d)$ and $\mathcal{N}_n^{k-}(d)$ to denote the index sets of discarded elements
\begin{subequations}
	\begin{align}
		\mathcal{N}_n^{k-}(d)&:=\underset{\mathcal{M}:\{\mathcal{M}\subset\mathcal{N}_n\cup\mathcal{B}_n, |\mathcal{M}|=q_n\}}{\arg\min}\sum_{m\in \mathcal{M}} [\bm{\theta}_m^k]_d, \\
		\mathcal{N}_n^{k+}(d)&:=\underset{\mathcal{M}:\{\mathcal{M}\subset\mathcal{N}_n\cup\mathcal{B}_n, |\mathcal{M}|=q_n\}}{\arg\max}\sum_{m\in \mathcal{M}} [\bm{\theta}_m^k]_d,
	\end{align}
\end{subequations}
where $[\bm{\theta}_m^k]_d$ represents the $d$-th element of $\bm{\theta}_m^k$. If we define $\mathcal{N}_n^{k*}(d):=\mathcal{N}_n-\mathcal{N}_n^{k+}-\mathcal{N}_n^{k-}$ as the set of the remaining indexes, the Byzantine-resilient update of dimension $d$ corresponding to  \eqref{equation:aggregate-decentralized-mean} is
\begin{equation}
	\label{equation:aggregate-decentralized-trimmedMean}
	[\bm{\theta}_n^{k+1}]_d = \frac{1}{N_n+{B_n}-2q_n+1}\!\!\sum_{m\in\mathcal{N}_n^{k*}(d)\cup \{n\}}\!\![\bm{\theta}_m^k]_d+\eta^k [\bm{g}_n^k(\bm{\theta}_n^k)]_d.
\end{equation}
The proposed \AlgorithmFullName{} is summarized in Algorithm \ref{algorithm:TD_TrimmedMean}.

\subsection{Why the trimmed mean works?}
\label{sec:effect-of-trimmed-mean}
We will discuss why trimmed mean is able to robustify the policy evaluation process. Intuitively, if the models sent by Byzantine agents are outliers, they will be discarded directly, so the output of trimmed mean will not go far away from the mean of true model. However, since trimmed mean aggregation discard $2q_n$ models at each coordinate and $q_n\ge B_n$, at least $q_n$ models from honest agents are discarded at the same time. 
Fortunately, trimmed mean can still approximate the true mean because information in the discarded models can be recover by that in $\mathcal{N}_n^{k*}(d)$.

Consider a simple scenario where a Byzantine agent $n'$ mixes into $\mathcal{N}_n^{k*}(d)$ and honest agents $\bm{\theta}_{n^{k+}}^k\in\mathcal{N}_n^{k+}(d)$ and $\bm{\theta}_{n^{k-}}^k\in\mathcal{N}_n^{k-}(d)$ are discarded. As a result, it satisfies that
\begin{align}
	[\bm{\theta}_{n^{k+}}^k]_d
	\le [\bm{\theta}_{n'}^k]_d
	\le [\bm{\theta}_{n^{k-}}^k]_d.
\end{align}
As a result, $[\bm{\theta}_{n'}^k]_d$ can be expressed as the linear combination of $[\bm{\theta}_{n^{k+}}^k]_d$ and $[\bm{\theta}_{n^{k-}}^k]_d$, that is (with a constant $y\in[0, 1]$)
\begin{align}
	[\bm{\theta}_{n'}^k]_d =
	y[\bm{\theta}_{n^{k+}}^k]_d
	+ (1-y)[\bm{\theta}_{n^{k-}}^k]_d.
\end{align}
This example shows that it is possible to recover the information in the discarded models from honest agents ($[\bm{\theta}_{n^{k+}}^k]_d$ and $[\bm{\theta}_{n^{k-}}^k]_d$) and the trimmed mean can be viewed as a linear combination of the models from neighbors.

While the trimmed mean has been used in other Byzantine-resilient algorithms \cite{blanchard2017MachineLW,yang2019bridge,guo2020towards}, most of them cannot give a finite-time upper bound of consensus and convergence rate. All of them require either the gradient of the objective function or its unbiased stochastic estimate. The most relevant works \cite{yang2019bridge} and \cite{vaidya2012matrix} are based on deterministic gradient descent. This makes our analysis of Byrd-TD($\lambda$) nontrivial.

\section{Theoretical Analysis}
\label{sec:analysis}

In this section, we first provide some important properties of TD($\lambda$), and then establish the convergence and consensus analysis for Byrd-TD($\lambda$). We will show the proof sketch in Section \ref{sec:proofsketch}, and leave the proofs to the appendices.

\subsection{Stationary point and asymptotic properties}
\label{sec:stationaryTD}


In the TD($\lambda$) recursion given by \eqref{equation:update-TD-normal}, a stationary point of $\{{\bm{\theta}}^k\}$ implies the increment to be $\bm{0}$. For the decentralized case, we characterize the stationary point first in terms of function value and then in terms of parameter. To do so, for any $\bm{\theta}\in\mathbb{R}^D$, define the global increment as
\begin{align}
	\bar{\bm{g}}^k(\bm{\theta}) &:=
	\frac{1}{N}\sum_{n\in\mathcal{N}}\bm{g}_n^k(\bm{\theta}) \nonumber\\
	&=\frac{1}{N}\sum_{n\in\mathcal{N}} \left(
	r_n^k+(\gamma\bm{\phi}(s^{k+1})-\bm{\phi}(s^{k}))^{\top}\bm{\theta}
	\right) \bm{z}^k \nonumber\\
	&=\mathbf{A}^k\bm{\theta}+\bar{\bm{b}}^k,
\end{align}
where ${\bm{g}}^k_n$ is given by \eqref{definition:TD_increment}, and the coefficients $\mathbf{A}^k$ and $\bar{\bm{b}}^k$ are given by
\begin{equation}
	\label{definition:Akbk}
	\mathbf{A}^k := \bm{z}^k\left(\gamma\bm{\phi}(s^{k+1})-\bm{\phi}(s^{k})\right)^{\top} \text{and}~~
	\bar{\bm{b}}^k := \frac{1}{N}\sum_{n\in\mathcal{N}} r_n^k\bm{z}^k.
\end{equation}

Similarly, we define $\bm{b}^k_n := r^k_n\bm{z}^k$ satisfying $\frac{1}{N}\sum_{n\in\mathcal{N}}\bm{b}^k_n=\bar{\bm{b}}^k_n$, so that the local update can be expressed as $\bm{g}_n^k(\bm{\theta})=\mathbf{A}^k\bm{\theta}+\bm{b}^k_n$.


Observe that $\mathbf{A}^k$ and $\bar{\bm{b}}^k$ are both determined by the state trajectory $\{s^0, s^1,\cdots, s^k\}$ generated from the Markov chain. Therefore, to characterize the stationarity of
$\{\bar{\bm{g}}^k(\bm{\theta})\}$,
we require the Markov chain to satisfy the following assumption.

\begin{Assumption}[Markov chain]
	\label{assumption:markovChain}
	The Markov chain is irreducible and aperiodic.
\end{Assumption}

As highlighted in \cite{tsitsiklis1997analysis}, Assumption \ref{assumption:markovChain} guarantees the existence of the stationary distribution $\rho$ of state $s^k$ and the expectation limits of $\mathbf{A}^k$ and $\bar{\bm{b}}^k$, as
\begin{equation}\label{eq:Ab}
	\lim_{k\to\infty}\mathbb{E}[\mathbf{A}^k]=\mathbf{A}^*
	\quad\text{and}\quad
	\lim_{k\to\infty}\mathbb{E}[\bar{\bm{b}}^k]=\bar{\bm{b}}^*.
\end{equation}
Therefore, the limit of the global increment can be defined as
\begin{equation}
	\lim_{k\to\infty}\mathbb{E}[\bar{\bm{g}}^k(\bm{\theta})] = \bar{\bm{g}}^*(\bm{\theta})~~~{\rm with}~~~	\bar{\bm{g}}^*(\bm{\theta}) := \mathbf{A}^*\bm{\theta}+\bar{\bm{b}}^*.
\end{equation}
It has been shown that $\mathbf{A}^*$ and $\bar{\bm{b}}^*$ have explicit expressions and $\mathbf{A}^*$ is a negative definite matrix \cite{tsitsiklis1997analysis}. In addition, $\mathbf{A}^*$ and $\bar{\bm{b}}^*$ will be attracted to their limits at a geometric mixing rate \cite{doan2019finite,bhandari2018finite}. These properties are presented in the supplementary, as well as the full version of this paper \cite{wu2020byzantine}. 

If there exists $\bm{\theta}^{\infty}_{\lambda}\in\mathbb{R}^D$ such that $\bar{\bm{g}}^*(\bm{\theta}^{\infty}_{\lambda})=\bm{0}$, we say that $\bm{\theta}^{\infty}_{\lambda}$ is a stationary point. We will show that our proposed Byrd-TD($\lambda$) converges to a neighborhood of $\bm{\theta}^{\infty}_{\lambda}$ due to Byzantine attacks. Note that although we can solve \eqref{equation:problem} with either gradient descent or TD($\lambda$), their stationary points are different. The stationary point of gradient descent minimizes $F(\bm{\theta})$, while in TD($\lambda$), different $\lambda$ leads to different $F(\bm{\theta}^{\infty}_{\lambda})$ because the function $\bar{\bm{g}}^*(\bm{\theta})$ depends on $\lambda$. This is one of the major challenges in analyzing Byrd-TD($\lambda$). It has been shown in \cite{tsitsiklis1997analysis} that $\bm{\theta}^{\infty}_{\lambda}$ satisfies
\begin{equation}
	\label{inequality:approximationAccuracy}
	\min_{\bm{\theta}} F(\bm{\theta})
	\le F(\bm{\theta}^{\infty}_{\lambda})
	\le \frac{1-\gamma\lambda}{1-\gamma}\min_{\bm{\theta}} F(\bm{\theta}).
\end{equation}
Larger $\lambda$ usually leads to smaller approximation error of $\bm{\theta}^{\infty}_{\lambda}$ in the original TD($\lambda$). However, as we can see later, this does not necessarily hold
when confronting with the Byzantine agents.



Denote the filtration containing all the information up to time $k$ by $\mathcal{F}^k=\{s^0, s^1,\cdots, s^k\}$. Unlike the analysis of stochastic optimization, the conditional expectation of $\bar{\bm{g}}^k(\bm{\theta})$ is usually not equal to $\bar{\bm{g}}^*(\bm{\theta})$ even when $k$ is large enough, because ${\cal P}(s^{k+1}=s|\mathcal{F}^k)$ is not equal to the stationary probability $\rho(s)$ even when $k$ goes to infinity. As a result, $\bar{\bm{g}}^k(\bm{\theta})$ is a biased estimator of $\bar{\bm{g}}^*(\bm{\theta})$.
This brings another of the major challenges to establish the convergence of Byrd-TD($\lambda$). To overcome this issue, we will resort to the geometric mixing time of Markov chain \cite{bremaud2013markov}.


\subsection{Main results}\label{sec:main-results}
Before presenting our main results, we first introduce several assumptions.
\begin{Assumption}[Bounded reward variation]
	\label{assumption:rewardVariation}
	For any transition from $s$ to $s'$, variation of the local reward at every honest agent with respect to the global reward is upper-bounded by
	\begin{equation}
		\label{equation:rewardVariation}
		\frac{1}{N}\sum_{n\in\mathcal{N}}\Big\|\mathcal{R}_n(s, s')-\frac{1}{N}\sum_{n'\in\mathcal{N}} \mathcal{R}_{n'}(s, s')\Big\|^2\le\delta^2.
	\end{equation}
\end{Assumption}

The quantity $\delta^2$ in Assumption \ref{assumption:rewardVariation} is the measurement of agent heterogeneity. When all agents have access to the global reward, $\delta^2$ is equal to $0$. As we will demonstrate in the analysis, a large $\delta^2$ increases the difficulty of defending against the Byzantine attacks.

\begin{Assumption}[Normalized features]
	\label{assumption:NormalizedFeature}
	Features are normalized such that $\|\bm{\phi}(s)\|\le 1,~\forall s\in\mathcal{S}$.
\end{Assumption}

Assumption \ref{assumption:NormalizedFeature} is standard in analyzing the TD-family algorithms \cite{doan2019finite,bhandari2018finite,srikant2019finite}.

We need an additional assumption on the network topology. Consider a set $\mathcal{H}_\mathcal{G}$ with cardinality $H_\mathcal{G}=|\mathcal{H}_\mathcal{G}|$ whose elements are the subgraphs of $\mathcal{G}$ obtained by removing all Byzantine agents with their edges, and removing any additional $q_n$ incoming edges at every honest node $n$.

\begin{Assumption}[Network connectivity]
	\label{assumption:networkStructure}
	For any subgraph $\mathcal{G}' \in \mathcal{H}_\mathcal{G}$, there exists at least one agent $n^*$ which has directed paths to all nodes in $\mathcal{G}'$. The length of paths is no more than $\tau_\mathcal{G}$.
\end{Assumption}

We call $n^*$ and $\tau_\mathcal{G}$ as \textit{source node} and \textit{network diameter}, respectively.
Examples of topologies that satisfy Assumption \ref{assumption:networkStructure} have been discussed in \cite{vaidya2012iterative}.

With these assumptions, we establish our main results. We first show all parameters in $\{\bm{\theta}^{k}_n\}$ will reach consensus albeit the Byzantine agents are biasing the learning process.
\begin{Theorem}
	\label{theorem:consensu-rate-theorem-global}
	Suppose Assumptions \ref{assumption:markovChain}--\ref{assumption:networkStructure} hold { and $q_n$ satisfies $B_n \le q_n < \frac{N_n}{3}$}. If we choose the decreasing step size $\eta^k$ in Algorithm \ref{algorithm:TD_TrimmedMean} appropriately satisfying
	\begin{subequations}
	\begin{align}
		& \label{condition:stepsize-consensus-1}
		16N^2\Big(\frac{1+\gamma}{1-\gamma\lambda}\Big)^2 (\eta^k)^2 
		\le \frac{1-\mu}{2} \\
		\label{condition:stepsize-consensus-2}
		&1\le
		\left(\frac{\eta^{k-1}}{\eta^{k}}\right)^2
		\le \frac{4}{3+\mu},
	\end{align}
	\end{subequations}
	and define $\bar{\bm{\theta}}^{k} := \frac{1}{N}\sum_{n\in\mathcal{N}} \bm{\theta}^{k}_n$, then the consensus error is
	\begin{equation}
		\label{inequality:consensu-rate-global}
		\frac{1}{N}\sum_{n\in\mathcal{N}}\left\|\bm{\theta}^{k}_n-\bar{\bm{\theta}}^{k}\right\|^2
		\le \frac{1}{2}C_1\mu_{\mathcal{G}}\left(\eta^{k}\right)^2,
	\end{equation}
	where the coefficients $\mu_{\mathcal{G}}$ and constant $C_1$ are defined as
	\begin{align}
		C_1 :=\frac{128N^2\delta^2}{(1-\gamma\lambda)^2},
		~~~~
		\mu_{\mathcal{G}} := \frac{1}{1-\mu},
	\end{align}
	and $\mu$ increases as the network diameter $\tau_\mathcal{G}$ increases.
	
\end{Theorem}

Theorem \ref{theorem:consensu-rate-theorem-global} implies that all honest agents eventually reach consensus, at the rate of $O\left((\eta^k)^2\right)$. A smaller network diameter $\tau_\mathcal{G}$ yields a faster convergence rate, as indicated by \eqref{inequality:consensu-rate-global}.

Next, we show that $\bm{\theta}^{k}_n$ on every honest agent $n$ converges to a neighborhood of $\bm{\theta}^{\infty}_\lambda$, the stationary point of the TD($\lambda$) recursion without Byzantine attacks.

\begin{Theorem}
	\label{theorem:convergence-trimmedMean}
	Suppose Assumptions \ref{assumption:markovChain}--\ref{assumption:networkStructure} hold and $q_n$ satisfies $B_n \le q_n < \frac{N_n}{3}$. If the step size $\eta^k$ in Algorithm \ref{algorithm:TD_TrimmedMean} satisfies $\eta^k=\frac{\eta}{k+k^0}$,
	where $k^0 > 0$ is a sufficiently large integer and $\eta > 0$, $\sigma_{\rm min}$ is the smallest singular value of $\mathbf{A}^*$, 
	then it holds
	\begin{align}
		\label{inequality:convergence_trimmedMean}
		\frac{1}{N}\!\sum_{n\in\mathcal{N}}\!\mathbb{E}\|&\bm{\theta}^{k+1}_n-\bm{\theta^{\infty}}_{\lambda}\|^2 \!\le\!
		C_1\mu_{\mathcal{G}}\left(\frac{\eta}{k+k^0}\right)^2
		+\frac{C_2}{(k+k^0)^\epsilon}
		+C_3 \varphi_\epsilon(k)
		+C_4 \frac{D_{\mathcal{G}}\delta^2}{(1-\gamma\lambda)^2},
	\end{align}
	where $\epsilon := \frac{\sigma_{\rm min}\eta}{2}$, $C_2$, $C_3$, and $C_4$ are positive constants, $D_{\mathcal{G}}$ is the degree of network unsaturation defined as
	\begin{equation}
		\label{definition:networkunsaturability}
		D_{\mathcal{G}} := \frac{N}{\min_{n\in\mathcal{N}}\{N_n+ {B_n}-2q_n+1\}}-1,
	\end{equation}
	and the function $\varphi_\epsilon(k)$ is defined as
	\begin{align}
		\label{definition:convergence-rate-dominant}
		\varphi_\epsilon(k) \!:=\!
		\left\{
		\begin{array}{r}
			\begin{aligned}
				&\frac{2\eta}{(\epsilon-1)^2}
				\frac{(\epsilon-1)\ln\big((k+1+k^0)/\eta\big)+1}{k+1+k^0},& &\!\!\!\!\!\epsilon>1,\\
				&\eta\frac{\ln\left((k+k^0)/\eta\right)^2}{k+k^0},& &\!\!\!\!\!\epsilon=1,\\
				& \frac{\eta\left((1-\epsilon)\ln\left(k^0/\eta\right)-1\right)}{(k^0-1)(1-\epsilon)^2}
				\left(\frac{k^0-1}{k+k^0}\right)^{\epsilon},& &\!\!\!\!\!0<\epsilon<1.
			\end{aligned}
		\end{array}
		\right.
	\end{align}
\end{Theorem}

Theorem \ref{theorem:convergence-trimmedMean} asserts that $\bm{\theta}^{k}_n$ on every honest agent $n$ converges to a neighborhood of $\bm{\theta}^\infty_\lambda$, with the rate dominated by the function $\varphi_\epsilon(k)$. When $\epsilon>1$, the algorithm converges at the rate of $O(\ln(k)/k)$, which is consistent with the rate in \cite{bhandari2018finite}. When $\epsilon=1$, the  algorithm converges at the rate of $O(\ln(k)^2/k)$, which is consistent with the rate in \cite{doan2019finite}. When $0<\epsilon<1$, the algorithm converges at the rate of $O(1/k^\epsilon)$. The radius of the neighborhood is bounded by a constant $O(D_{\mathcal{G}}\delta^2)$. To the best of our knowledge, this tight analysis is new.

In addition to the learning error caused by Byzantine attacks as shown in \eqref{inequality:convergence_trimmedMean}, the following corollary gives the overall asymptotic learning error.

\begin{Corollary}
	\label{corollary:asymptoticTDerror}
	Under the same condition as that in Theorem \ref{theorem:convergence-trimmedMean}, the asymptotic learning error of Algorithm \ref{algorithm:TD_TrimmedMean} is bounded by
	\begin{equation}
		\label{inequality:asymptoticTDerror}
		\lim_{k\to\infty}\sup\frac{1}{N}\sum_{n\in\mathcal{N}} F(\bm{\theta}^{k+1}_n)
		\le
		2C_4 \frac{D_{\mathcal{G}}\delta^2}{(1-\gamma\lambda)^2}
		+2\frac{1-\gamma\lambda}{1-\gamma}\min_{\bm{\theta}}F(\bm{\theta}).
	\end{equation}
\end{Corollary}

The first term at the right-hand side of \eqref{inequality:asymptoticTDerror} comes from the Byzantine attacks and the use of trimmed mean, while the second term is due to the intrinsic error of linear approximation. Tuning the parameter $\lambda$ helps achieve the best accuracy.


\begin{Remark}[Effect of network topology]  
	Two important parameters that determine the convergence rate and the asymptotic learning error are the network diameter $\tau_\mathcal{G}$ and the degree of network unsaturation $D_{\mathcal{G}}$. 
	The dependence on the network diameter $\tau_\mathcal{G}$ has been shown in Theorem \ref{theorem:consensu-rate-theorem-global}. 
	For $D_{\mathcal{G}}$ in \eqref{definition:networkunsaturability}, we observe that the denominator of its first term is $\min_{n\in\mathcal{N}}\{N_n+ {B_n}-2q_n+1\}$, which corresponds to the bottleneck of the network with the smallest number of honest in-neighbors. 
	Ideally, every honest agent $n$ will choose $q_n$ large enough to defend against Byzantine neighbors ($q_n \geq B_n$), but not too large such that $D_{\mathcal{G}}$ is sufficiently small. When the number of Byzantine agents $B$ is finite, the number of honest agents $N$ goes to infinity, and the network is complete, the degree of saturation is sufficiently high. In this case,  $D_{\mathcal{G}}$ goes to $0$, and so as the asymptotic learning error. On the contrary, when the network is sparse, $\min_{n\in\mathcal{N}}\{N_n+ {B_n}-2q_n+1\}$ becomes small such that $D_{\mathcal{G}}$ and the asymptotic learning error are both $O(N)$.
	In short, dense networks are easier for learning the value function in terms of both consensus and learning error. However, more communication links lead to more overhead, leaving a trade-off for network designers.
\end{Remark}

\begin{Remark}[Selection of parameter $\lambda$]
	According to the original TD($\lambda$) analysis, $\lambda$ balances the approximation accuracy and the convergence rate. As demonstrated in \eqref{inequality:approximationAccuracy}, $\lambda=1$ leads to the best approximation value, while $\lambda=0$ leads to the fastest convergence rate. Due to the existence of Byzantine attacks, the influence of $\lambda$ becomes more complicated in the proposed Byrd-TD($\lambda$). At the right-hand side of \eqref{inequality:asymptoticTDerror}, when $\lambda$ increases from $0$ to $1$, the first term increases but the second term decreases.
	Recall that $\delta^2$ is the measurement of reward heterogeneity and $D_{\mathcal{G}}$ is the degree of network unsaturation. We suggest to choose a small $\lambda$ if the reward heterogeneity is high or the network is unsaturated.
\end{Remark}

\section{Proof Sketch of Main Results}
\label{sec:proofsketch}
In this section, we will give the key idea of proving the main results. As explained in Section \ref{sec:effect-of-trimmed-mean}, we can view the trimmed mean as a linear combination of models, which enables us to express the trimmed mean aggregation in a compact form.

To see so, we first introduce some notations. Concatenating all $\bm{\theta}_n^k$ and $\bm{g}^k_n(\bm{\theta}_n^k)$, we  define $\mathbf{\Theta}^k$, $\mathbf{G}^k(\mathbf{\Theta}^k)$ respectively as
\begin{align}
	\label{definition:matrix-Theta}
	\mathbf{\Theta}^k &:= \left[
	\bm{\theta}_1^k, \cdots,
	\bm{\theta}_n^k, \cdots, \bm{\theta}_{N}^k
	\right]^{\top}\in\mathbb{R}^{N\times D} \nonumber \\
	\mathbf{G}^k(\mathbf{\Theta}^k) &:= \left[
	\bm{g}^k_1(\bm{\theta}_1^k), \cdots,
	\bm{g}^k_n(\bm{\theta}_n^k), \cdots, \bm{g}^k_N(\bm{\theta}_N^k)
	\right]^{\top}\!\!\in\mathbb{R}^{N\times D}.
\end{align}


In addition, the observation in Section \ref{sec:effect-of-trimmed-mean} implies the trimmed mean is a coordinate-wise linear combination of the model from several neighbors, so it is natural to express the update \eqref{equation:aggregate-decentralized-trimmedMean} in a compact form, given by
\begin{align}
	\label{equation:update-trimmedMean-matrix}
	\mathbf{\Theta}^{k+1}_d=\mathbf{Y}^k(d)\mathbf{\Theta}^{k}_d+\eta^k \mathbf{G}^k_d(\mathbf{\Theta}^k),
\end{align}
where $\mathbf{\Theta}^k_d$ and $\mathbf{G}^k_d$ represent the $d$-th columns of $\mathbf{\Theta}^k$ and $\mathbf{G}^k$, respectively, and $\mathbf{Y}^k(d)\in\mathbb{R}^{N\times N}$ is a proper matrix, such that $[\mathbf{Y}^k(d)]_{nm}> 0$ if and only if $(n,m)\in\mathcal{E}^k(d)$. In section \ref{sec:contructionY}, we will construct $\mathbf{Y}^k(d)$ explicitly.

\subsection{Decentralized TD($\lambda$) over a varying directed network}
\label{sec:time-varying-main-result}
With the trimmed mean aggregation, at each iteration $k$, every agent discards some messages that may come from different agents. As a result, the communication graph is essentially a directed and time-varying subgraph of $\mathcal{G}$ (even $\mathcal{G}$ is undirected).
Therefore, we analyze Byrd-TD($\lambda$) by analyzing a more general algorithm --- decentralized TD($\lambda$) over a time-varying directed network, which by itself is new, and of independent interest. 

By generalizing \cite[Claim 2]{vaidya2012matrix}, under conditions described below, we will show the consensus and convergence of TD($\lambda$) over a time-varying directed network. 
Before describing these properties, we introduce some definitions.
\begin{Definition}
	(Stochastic vector and matrix) A vector $\bm{x}$ is termed a \textit{stochastic vector} if all its elements are in the range $[0, 1]$ and sum up to $1$.
	A matrix is termed a \textit{row-stochastic matrix} if all its rows are stochastic vectors.
\end{Definition}

Then, the conditions $\mathbf{Y}^k(d)$ should satisfy include
\begin{enumerate}
	\item[(c1)] $\mathbf{Y}^k(d)$ is a row-stochastic matrix.
	\item[(c2)] $[\mathbf{Y}^k(d)]_{nn} = \frac{1}{N_n+B_n-2q_n+1}$ for all $n\in\mathcal{N}$.
	\item[(c3)] $[\mathbf{Y}^k(d)]_{nm}\neq 0$ only if $(m, n)\in\mathcal{E}$ or $n=m$.
	\item[(c4)] At least $N_n-q_n+1$ elements in the $n$-th row of $\mathbf{Y}^k(d)$ are lower bounded by a positive constant, written as $\mu_0$.
	\item[(c5)] All elements in $\mathbf{Y}^k(d)$ are upper bounded by
	$\frac{1}{\min_{n\in\mathcal{N}}\{N_n+B_n-2q_n+1\}}$.
	\item[(c6)] There exists one column in $(\mathbf{Y}^k(d))^{\tau_\mathcal{G}}$ such that all elements in this column are non-zero and the there exist no more than $H_\mathcal{G}$ types of $\mathbf{Y}^k(d)$.
\end{enumerate}

With (c1)--(c6), we can derive the consensus rate of TD($\lambda$) over a time-varying directed network.
\begin{Theorem}
	\label{theorem:consensu-rate-theorem-global-time-varying}
	Suppose Assumptions \ref{assumption:markovChain}--\ref{assumption:NormalizedFeature}, Conditions (c1)--(c6) of $\mathbf{Y}^k(d)$ hold and $q_n$ satisfies $B_n \le q_n < \frac{N_n}{3}$ . If we choose the step size $\eta^k$ in Algorithm \ref{algorithm:TD_TrimmedMean} satisfying condition \eqref{condition:stepsize-consensus-1} and \eqref{condition:stepsize-consensus-2}, then the consensus error at iteration $k$ satisfies
	\begin{equation}
		\label{inequality:consensu-rate-global-time-varying}
		\frac{1}{N}\sum_{n\in\mathcal{N}}\left\|\bm{\theta}^{k}_n-\bar{\bm{\theta}}^{k}\right\|^2
		\le \frac{1}{2}C_1\mu_{\mathcal{G}}\left(\eta^{k}\right)^2.
	\end{equation}
\end{Theorem}

Theorem \ref{theorem:consensu-rate-theorem-global-time-varying} implies that all agents eventually reach consensus, at the rate of $O\left((\eta^k)^2\right)$.

Next, we show that $\bm{\theta}^{k}_n$ on every agent $n$ converges to a neighborhood of the fixed point $\bm{\theta}^{\infty}_\lambda$.

\begin{Theorem}
	\label{theorem:convergence-trimmedMean-time-varying}
	Suppose Assumptions \ref{assumption:markovChain}--\ref{assumption:NormalizedFeature}, Conditions (c1)--(c6) of $\mathbf{Y}^k(d)$ hold and $q_n$ satisfies $B_n \le q_n < \frac{N_n}{3}$ . If the step size $\eta^k$ in Algorithm \ref{algorithm:TD_TrimmedMean} satisfies $\eta^k=\frac{\eta}{k+k^0}$,
	where $k^0 > 0$ is a sufficiently large number such that step size $\eta^k$ satisfies the condition \eqref{condition:stepsize-consensus-1}, \eqref{condition:stepsize-consensus-2} and $\eta > 0$,
	then it holds
	\begin{align}
		\label{inequality:convergence_trimmedMean-time-varying}
		\frac{1}{N}\sum_{n\in\mathcal{N}}\mathbb{E}&\|\bm{\theta}^{k+1}_n-\bm{\theta}^{\infty}_{\lambda}\|^2 \le
		C_1\mu_{\mathcal{G}}\left(\frac{\eta}{k+k^0}\right)^2
		+\frac{C_2}{(k+k^0)^\epsilon}
		\!+C_3 \varphi_\epsilon(k)
		+C_4 \frac{D_{\mathcal{G}}\delta^2}{(1-\gamma\lambda)^2},
	\end{align}
	where $D_{\mathcal{G}}$ and the function $\varphi_\epsilon(k)$ have been defined in \eqref{definition:networkunsaturability} and \eqref{definition:convergence-rate-dominant}, respectively.
\end{Theorem}

Theorems \ref{theorem:consensu-rate-theorem-global-time-varying} and \ref{theorem:convergence-trimmedMean-time-varying} are the cornerstones to prove Theorems \ref{theorem:consensu-rate-theorem-global} and \ref{theorem:convergence-trimmedMean}, respectively. The remaining task is to show the effect of trimmed mean in Byrd-TD($\lambda$) can be described by a well-defined matrix $\mathbf{Y}^k(d)$ satisfying (c1)--(c6).

\subsection{Construction of $\mathbf{Y}^k(d)$}
\label{sec:contructionY}
By extending \cite{vaidya2012matrix}, we construct $\mathbf{Y}^k(d)\in\mathbb{R}^{N\times N}$ as follows. We use $B^{k*}_n(d) := |\mathcal{N}_n^{k*}(d)\cap\mathcal{B}|$ to denote the number of Byzantine agents in $\mathcal{N}_n^{k*}(d)$ and define $N^*_n := N_n+B_n-2q_n+1$ for simplicity.

We construct $\mathbf{Y}^k(d)$ under two cases.

\noindent
\textbf{Case 1: agent $n$ estimates $q_n$ correctly.}

This means $q_n-B_n+B^{k*}_n(d)=0$.
Since both $q_n-B_n$ and $B^{k*}_n(d)$ are larger than or equal to $0$, the condition $q_n-B_n+B^{k*}_n(d)=0$ implies $q_n-B_n=0$ and $B^{k*}_n(d)=0$, meaning that the honest agent $n$ can estimate the number of Byzantine neighbors correctly and eliminate all Byzantine neighbors. For any of honest neighbor $m$, the weight should be chosen to $\frac{1}{N_n+B_n-2q_n+1}=\frac{1}{N_n^*}$. 

Hence, $\mathbf{Y}^k(d)$ can be constructed by
\begin{align}
	[\mathbf{Y}^k(d)]_{nm} = \begin{cases}
		\frac{1}{N^*_n},~~~&{\rm if}~m=n~{\rm or}~m\in \mathcal{N}_n^{k*}(d)\cap\mathcal{N}_n\\
		~~  0,~~~&{\rm otherwise}.
	\end{cases}
\end{align}

The main different between our construction and that in \cite{vaidya2012matrix} is how to choose the weights in Case 2.

\noindent\textbf{Case 2: agent $n$ estimates $q_n$ incorrectly.}

This means that $q_n\!-\!B_n+B^{k*}_n(d)>0$.
When agent $n$ overestimates the number of Byzantine neighbors, we construct $\mathbf{Y}^k(d)$ tactfully. We can verify that\footnote{As a matter of fact, there are at least $q_n-B_n+B^{k*}_n(d)$ and at most $q_n$ honest agents in both $\mathcal{N}_n^{k+}(d)$ and $\mathcal{N}_n^{k-}(d)$.}
at least $q_n-B_n+B^{k*}_n(d)$ honest agents are contained in both $\mathcal{N}_n^{k+}(d)$ and $\mathcal{N}_n^{k-}(d)$.
Therefore, we define subsets $\mathcal{L}_n^{k+}(d)\subseteq\mathcal{N}_n^{k+}(d)$ and  $\mathcal{L}_n^{k-}(d)\subseteq\mathcal{N}_n^{k-}(d)$ where $|\mathcal{L}_n^{k+}(d)|=|\mathcal{L}_n^{k-}(d)|=q_n-B_n+B^{k*}_n(d)$.



When $\mathcal{N}_n^{k*}(d)\cap\mathcal{B}_n\neq\emptyset$ (then $B^{k*}_n(d) \neq 0$), for any $n^{k+}\in\mathcal{L}_n^{k+}(d)$ and $n^{k-}\in\mathcal{L}_n^{k-}(d)$, for any $n'\in\mathcal{N}_n^{k*}(d)\cap\mathcal{B}_n$, the element $[\bm{\theta}_{n'}^k]_d$ satisfies
\begin{align}
	[\bm{\theta}_{n^{k+}}^k]_d
	\le [\bm{\theta}_{n'}^k]_d
	\le [\bm{\theta}_{n^{k-}}^k]_d.
\end{align}
As a result, $[\bm{\theta}_{n'}^k]_d$ can be expressed as
\begin{align}
	[\bm{\theta}_{n'}^k]_d =
	y(n',n^{k+})[\bm{\theta}_{n^{k+}}^k]_d
	+ y(n',n^{k-})[\bm{\theta}_{n^{k-}}^k]_d,
\end{align}
where $y(n',n^{k+})$ and $y(n',n^{k-})$ are constants satisfying $0\le y(n',n^{k+}), y(n',n^{k-}) \le 1$ and $y(n',n^{k+})+y(n',n^{k-})=1$. In fact, $y(n',n^{k+})$ and $y(n',n^{k-})$ both depend on $k$ and $d$, but we ignore the subscripts for simplicity.

Similarly, defining the constant 
\begin{align}
	\label{definition:ckn}
	c_{n}^k(d) := \frac{q_n-B_n}{N_n+B_n-2q_n-B^{k*}_n(d)} < 1,
\end{align}
we can decompose $[\bm{\theta}_{n'}^k]_d$ for $n'\in\mathcal{N}_n^{k*}(d)\cap\mathcal{N}_n$ by
\begin{align}
	&	[\bm{\theta}_{n'}^k]_d 
	=\left(
	1-c_{n}^k(d)
	\right) [\bm{\theta}_{n'}^k]_d
	+ c_{n}^k(d)
	\left(y(n',n^{k+})[\bm{\theta}_{n^{k+}}^k]_d
	+ y(n',n^{k-})[\bm{\theta}_{n^{k-}}^k]_d\right).
\end{align}
Inspired by this, we can construct $\mathbf{Y}^k(d)$ as
\begin{align}
	\label{eq.Y2}
	[\mathbf{Y}^k(d)]_{nm} =
	\begin{cases}
		\frac{1}{N^*_n},~~~~&{\rm if}~m=n\\
		\left(1-c_{n}^k(d) \right)\frac{1}{N^*_n},~~~~&{\rm if}~m\in \mathcal{N}_n^{k*}(d)\cap\mathcal{N}_n\\
		\frac{\sum_{n'\in\mathcal{N}_n^{k*}(d)\cap\mathcal{N}_n} 
			y(n',m)\frac{1}{N^*_n} c_{n}^k(d)}{q_n-B_n+B^{k*}_n(d)} \\
		~+\frac{\sum_{n'\in\mathcal{N}_n^{k*}(d)\cap\mathcal{B}_n}    y(n',m)\frac{1}{N^*_n}}{q_n-B_n+B^{k*}_n(d)},&{\rm if}~m\in\mathcal{L}_n^{k+}(d)
		\cup\mathcal{L}_n^{k-}(d) \\
		0,~~~~&{\rm otherwise}.
	\end{cases}
\end{align}


With the definition of $\mathbf{Y}^k(d)$ in \eqref{eq.Y2}, we can write Byrd-TD($\lambda$) update in the form of \eqref{equation:update-trimmedMean-matrix}, where $\mathbf{Y}^k(d)$ satisfies conditions (c1)--(c6).

\section{Numerical Experiments}
We test the proposed Byrd-TD($\lambda$) on the \emph{cooperative navigation} task that modifies the one in \citep{lowe2017nips}. The goal of every agent is to cover its target landmark. An agent's local reward is computed by calculating the distance between itself and its target landmark and, is further penalized if it collides with other agents. In every step, an agent chooses an action among $\{up ,left, right, down\}$ according to its policy. The local reward of each agent depends on the distance between the agent and its goal landmark, and will be penalized if the agent collides with others. Each agent does not know others' landmarks, and thus the local reward functions are kept private. Unlike the task in \citep{lowe2017nips}, where every agent does not have a specific target landmark but has access to all local rewards, our setting is fully decentralized. The policy to be evaluated is uniform random. The performance metrics are the mean squared Bellman error (\texttt{MSBE}) and the mean consensus error (\texttt{MCE}). At every step, every honest agent has a local squared Bellman error (\texttt{SBE}) defined as
\begin{align}
    \label{definition:SBE}
    \texttt{SBE} \left(\{\bm{\theta}_n^k\}_{n=1}^{\mathcal{N}}, s^k \right) := \frac{1}{N}\sum_{n\in\mathcal{N}}
    \left(
    V_{\bm{\theta}_n^k}(s^k) - \sum_{s'\in\mathcal{S}} 
    \mathcal{P}(s^k, s')\big(
    \mathcal{R}_n(s^k, s') + \gamma V_{\bm{\theta}_n^k}(s')\big)
    \right)^2,
\end{align}
and a local consensus error (\texttt{CE}) defined as
\begin{align}
    \label{definition:CE}
    \texttt{CE} \left(\{\bm{\theta}_n^k\}_{n=1}^{\mathcal{N}} \right) := \frac{1}{N}\sum_{n\in\mathcal{N}}\|\bm{\theta}^{k}_n-\bar{\bm{\theta}}^{k}\|^2.
\end{align}

\texttt{MSBE} and \texttt{MCE} are calculated by averaging the local \texttt{SBE}s and \texttt{CE}s over all honest agents and all previous steps, respectively, given by
\begin{align}
    \texttt{MSBE} := \frac{1}{k}\sum_{\kappa=1}^{k} \texttt{SBE} \left(\{\bm{\theta}_n^\kappa\}_{n=1}^{\mathcal{N}}, s^\kappa \right)~~~~~~{\rm and}~~~~~~
    \texttt{MCE} := \frac{1}{k}\sum_{\kappa=1}^{k} \texttt{CE} \left(\{\bm{\theta}_n^\kappa\}_{n=1}^{\mathcal{N}} \right)
\end{align}

\begin{figure*}[t]
		\hspace{-0.55cm}
		\includegraphics[width=1.05\textwidth]{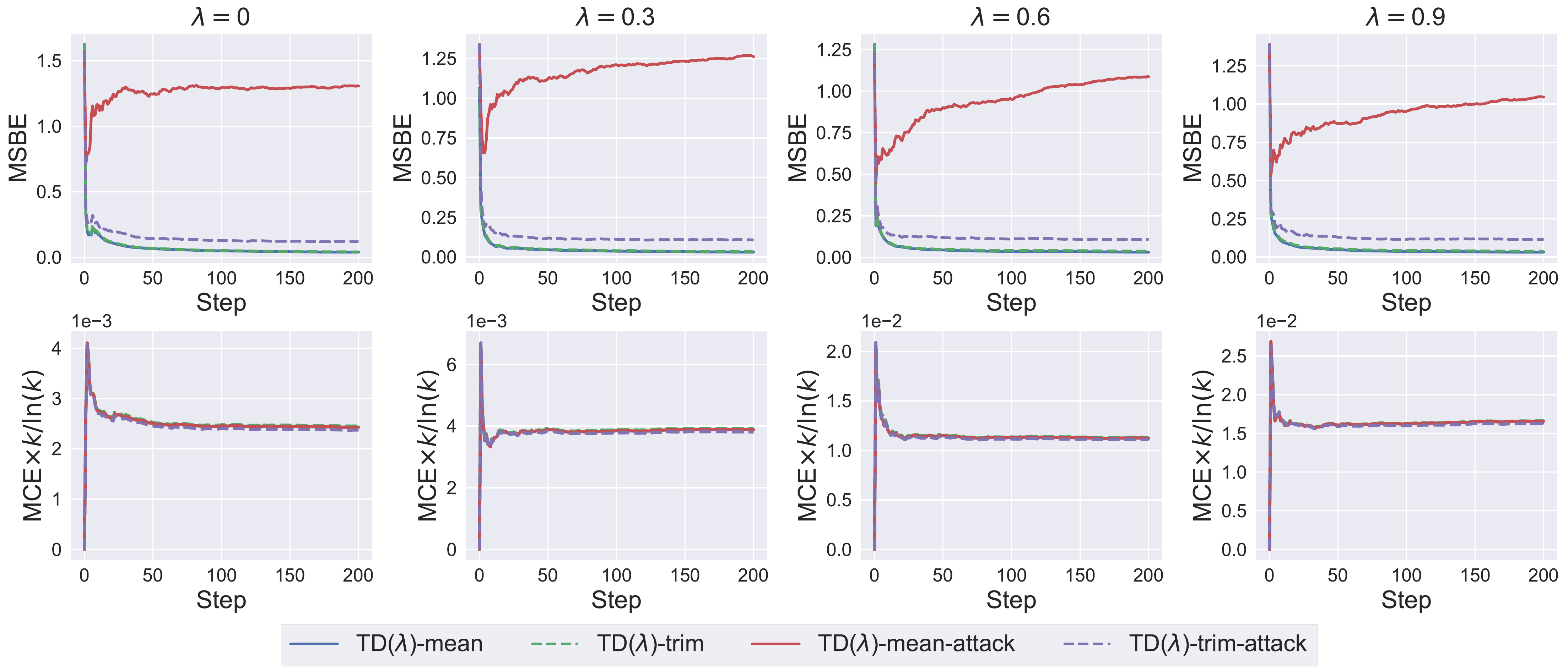}
	\vspace*{-0.3cm}
	\caption{\texttt{MSBE} and \texttt{MCE}$\times k/\ln(k)$ versus step $k$ under different $\lambda$ in a complete network.}
	\label{fig:complete}
	\vspace*{-0.1cm}
\end{figure*}

\begin{figure*}[t]
		\hspace{-0.55cm}
		\includegraphics[width=1.05\textwidth]{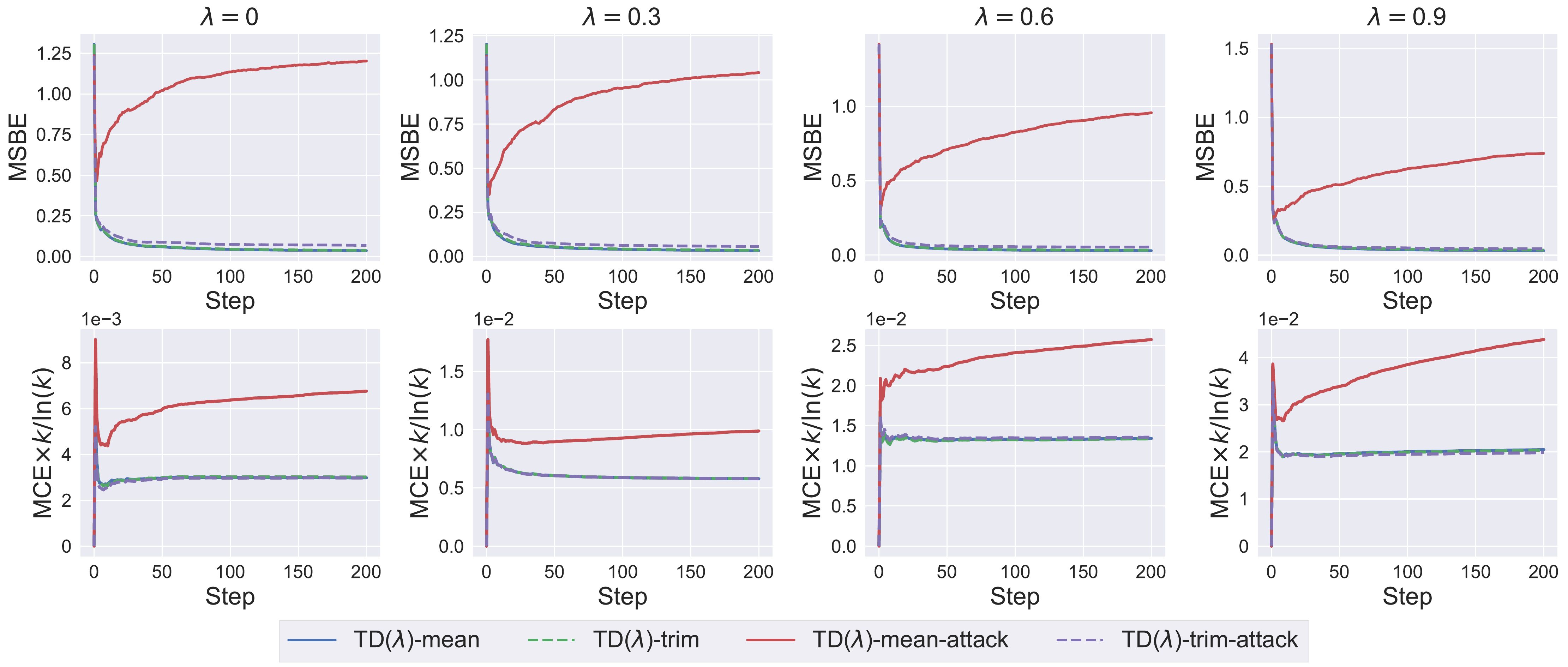}
	\vspace*{-0.3cm}
	\caption{\texttt{MSBE} and \texttt{MCE}$\times k/\ln(k)$ versus step $k$ under different $\lambda$ in an Erdos-Renyi network.}
	\label{fig:renyi}
	\vspace*{-0.5cm}
\end{figure*}

\textbf{Byzantine-resilience.} We compare four algorithms: \textit{mean} and \textit{trim} refer to that the honest agents update their parameters via mean and trimmed mean using their neighboring parameters without Byzantine attacks, respectively; \textit{mean-attack} and \textit{trim-attack} refer to those under Byzantine attacks. That is, \textit{mean} and \textit{mean-attack} correspond to the original decentralized TD($\lambda$), while \textit{trim} and \textit{trim-attack} correspond to  Byrd-TD($\lambda$). 
 {We consider sign flipping attacks here.}
We examine the impact of different $\lambda$. For $\lambda=0, 0.3, 0.6$, the step size is $\eta^k=\frac{0.1}{\sqrt{k}}$. For $\lambda=0.9$, the step size is $\eta^k=\frac{0.05}{\sqrt{k}}$. The first underlying network is complete, with $7$ honest and $2$ Byzantine agents, and $q_n = 2$. The second is an Erdos-Renyi network with $9$ agents. Every pair of agents are neighbors with probability $0.7$, and every agent is Byzantine with probability $0.2$. We set $q_n$ as the number of Byzantine agents. A total of 10 random Erdos-Renyi graphs are generated to calculate the averaged \texttt{MSBE} and \texttt{MCE}.  The Byzantine agents are normal in \textit{mean} and \textit{trim}, but adopt sign flipping attacks in \textit{mean-attack} and \textit{trim-attack}. Here, \textit{sign flipping} attacks mean that Byzantine agents send negative values of their true parameters to its neighbors. 

As shown in Figures \ref{fig:complete} and \ref{fig:renyi}, when the Byzantine agents are absent, decentralized TD($\lambda$) with mean and trimmed mean aggregation rules both work well. At presence of Byzantine attacks, decentralized TD($\lambda$) with mean aggregation fails, but the proposed trimmed mean aggregation rule is still robust. With particular note, as predicted by Theorem \ref{theorem:consensu-rate-theorem-global}, \texttt{CE} is in the order of $O((\eta^k)^2) = O(\frac{1}{k})$, such that \texttt{MCE} is in the order of $O(\frac{1}{k}\sum_{\kappa=1}^{k}(\eta^\kappa)^2) = O(\ln(k)/k)$ and \texttt{MCE}$\times k/\ln(k)$ becomes a horizontal line asymptotically. Our numerical experiments corroborate this theoretical result.

\begin{wrapfigure}{r}{.5\linewidth}
    \vspace{-0.4cm}
    \def\epsfsize#1#2{0.4#1}
    \includegraphics[width=1\linewidth]{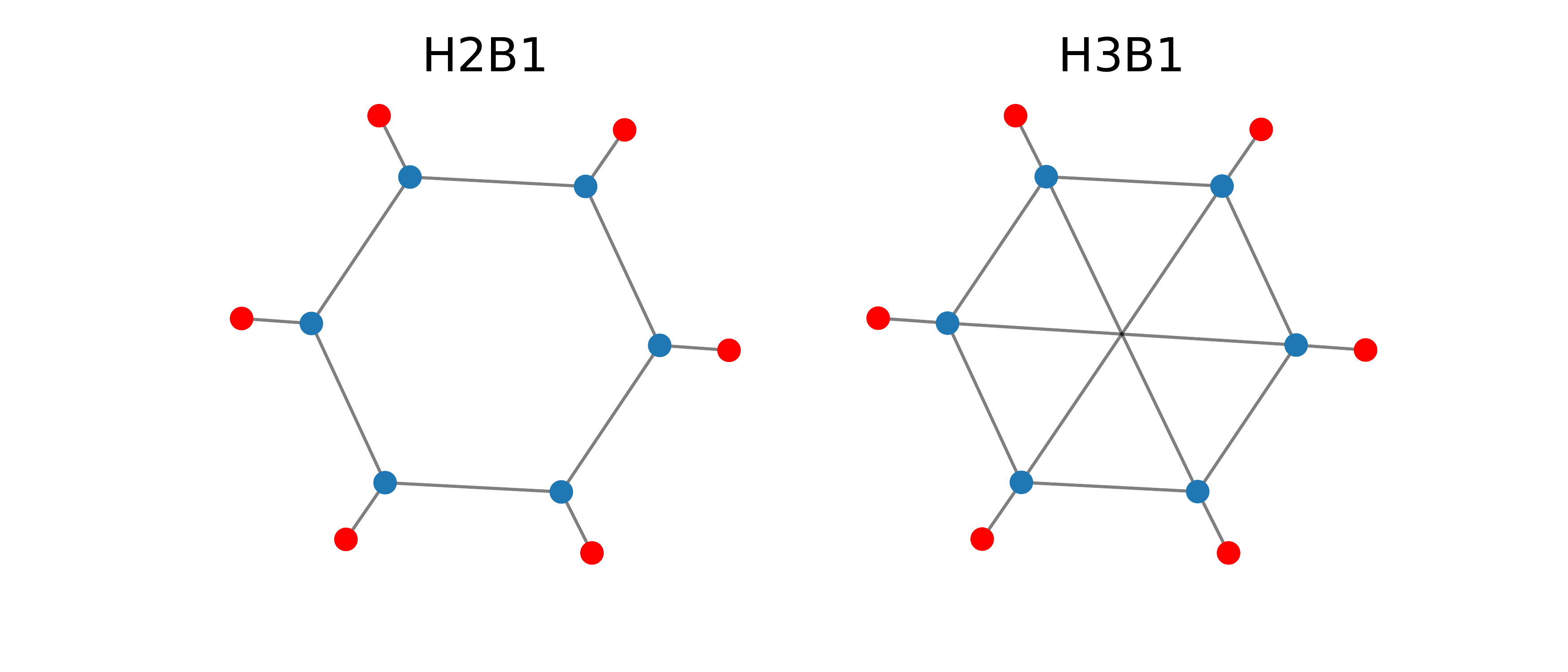}
	\vspace*{-0.75cm}
	\caption{Each honest agent has one Byzantine neighbor.}
    \label{fig:oneByzantine}
    \vspace*{-1.0cm}
\end{wrapfigure}

\begin{figure*}[t]
		\hspace{-0.55cm}
		\includegraphics[width=1.05\textwidth]{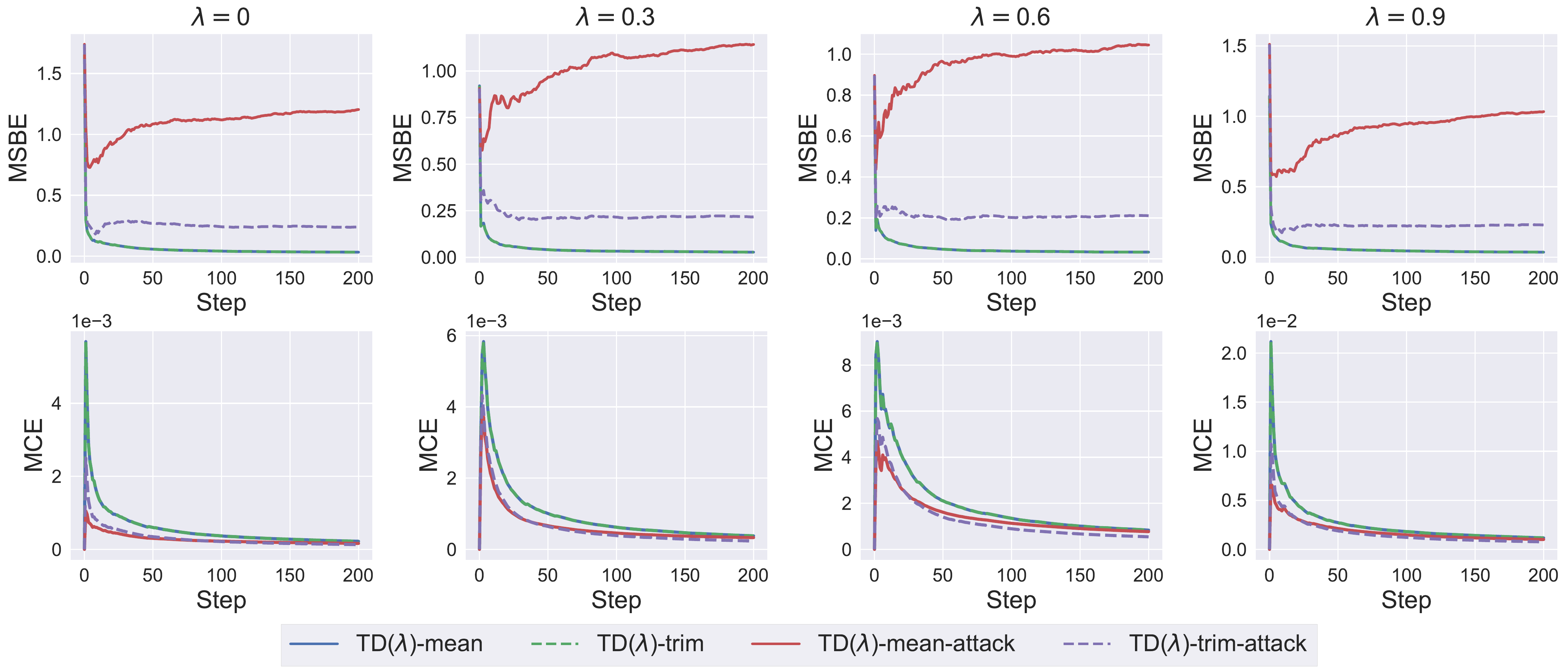}
	\vspace*{-0.3cm}
	\caption{\texttt{MSBE} and \texttt{MCE} under sign flipping attacks in H2B1 network shown in Figure \ref{fig:oneByzantine}.}
	\label{fig:h2b1}
	\vspace*{-0.3cm}
\end{figure*}

\begin{figure*}[t]
		\hspace{-0.55cm}
		\includegraphics[width=1.05\textwidth]{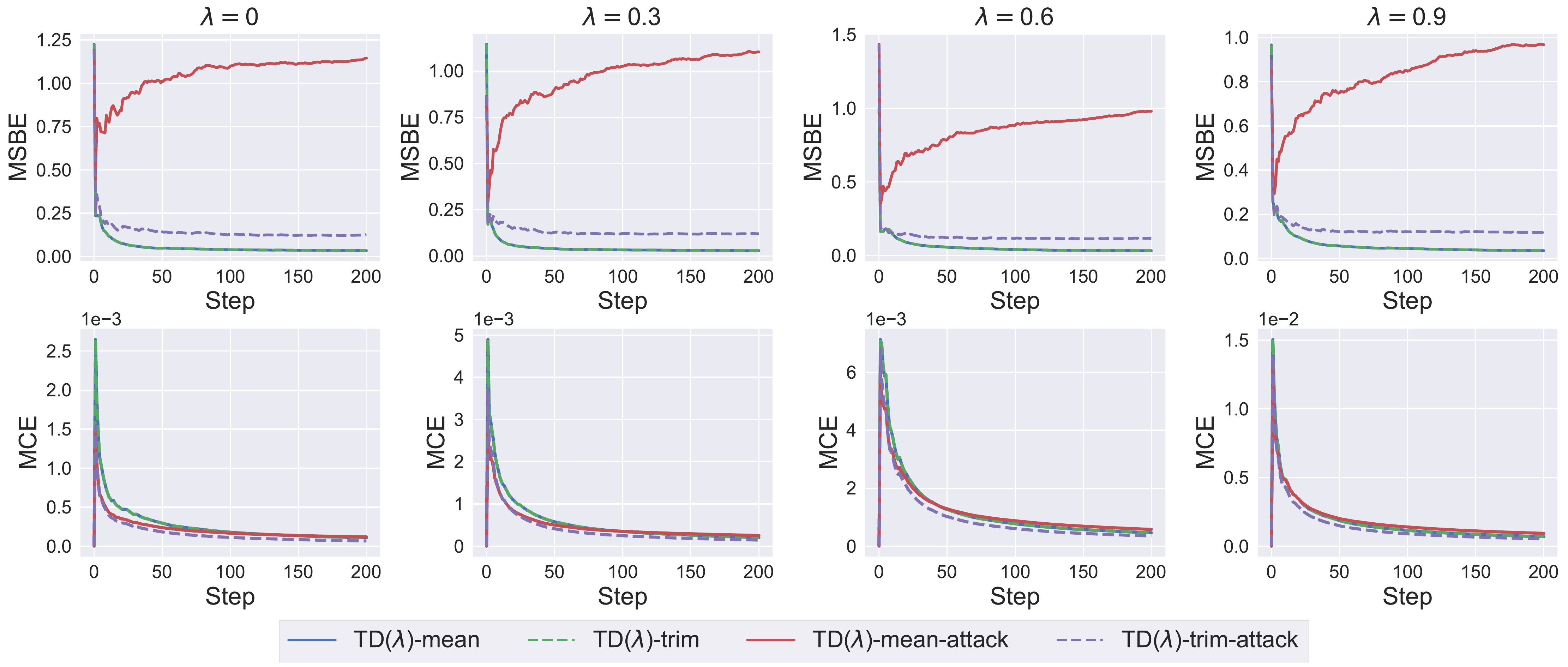}
	\vspace*{-0.3cm}
	\caption{\texttt{MSBE} and \texttt{MCE} under sign flipping attacks in H3B1 network shown in Figure \ref{fig:oneByzantine}.}
	\label{fig:h3b1}
	\vspace*{-0.3cm}
\end{figure*}

\textbf{Network topology.} We test the algorithms over different topologies. We have presented the performance in complete and Erdos-Renyi graphs under sign flipping attacks in Figures \ref{fig:complete} and \ref{fig:renyi}. 
In Figures \ref{fig:h2b1} and \ref{fig:h3b1}, we test the algorithm in several decentralized networks shown in Figure \ref{fig:oneByzantine}. In these experiments, $q_n=1$. For $\lambda=0, 0.3, 0.6$, the step size is set as $\eta^k=\frac{0.1}{\sqrt{k}}$. For $\lambda=0.9$, $\eta^k=\frac{0.05}{\sqrt{k}}$.
Byzantine agents are normal in \textit{mean} and \textit{trim}, but adopt sign flipping attacks in \textit{mean-attack} and \textit{trim-attack}.

\begin{figure*}[t]
		\hspace{-0.55cm}
		\includegraphics[width=1.05\textwidth]{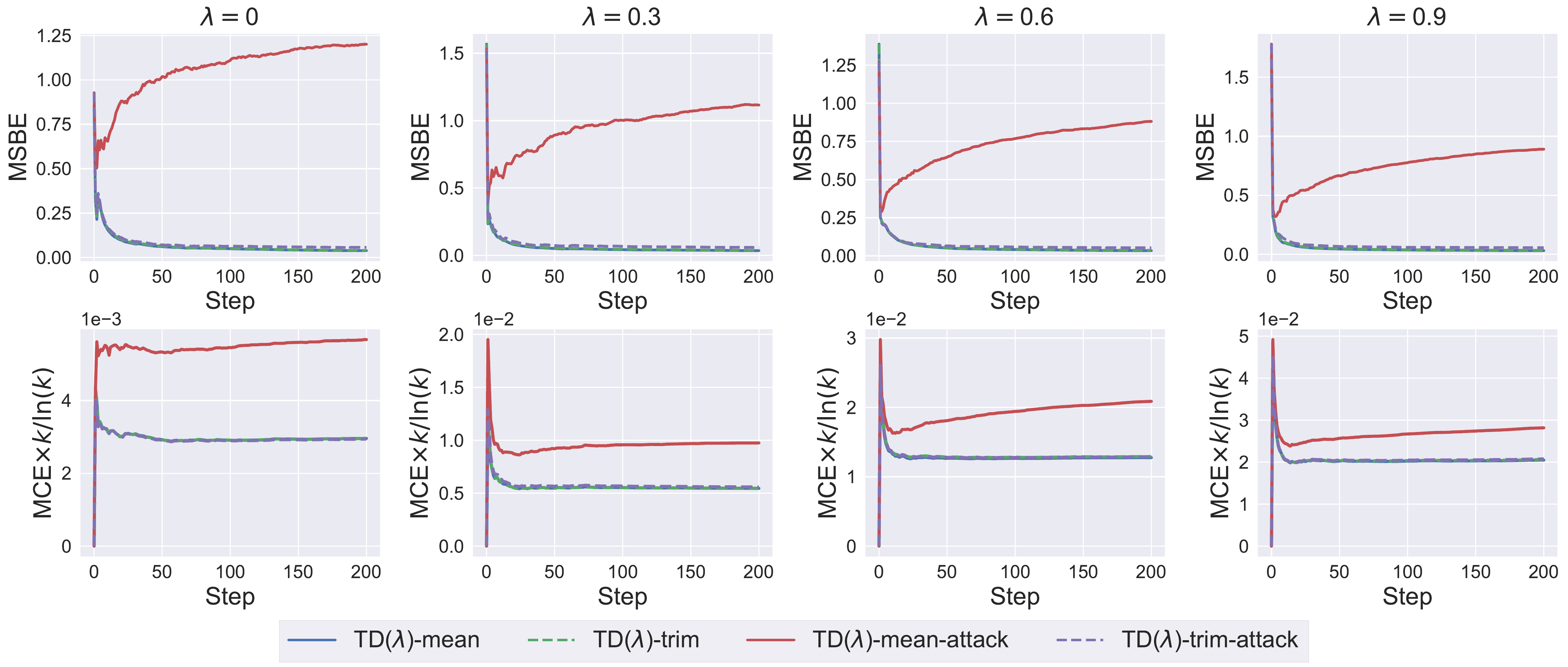}
	\vspace*{-0.3cm}
	\caption{\texttt{MSBE} and \texttt{MCE} under same value attacks in Erdos-Renyi graph.}
	\label{fig:renyi-a0}
	\vspace*{-0.3cm}
\end{figure*}

\begin{figure*}[t]
		\hspace{-0.55cm}
		\includegraphics[width=1.05\textwidth]{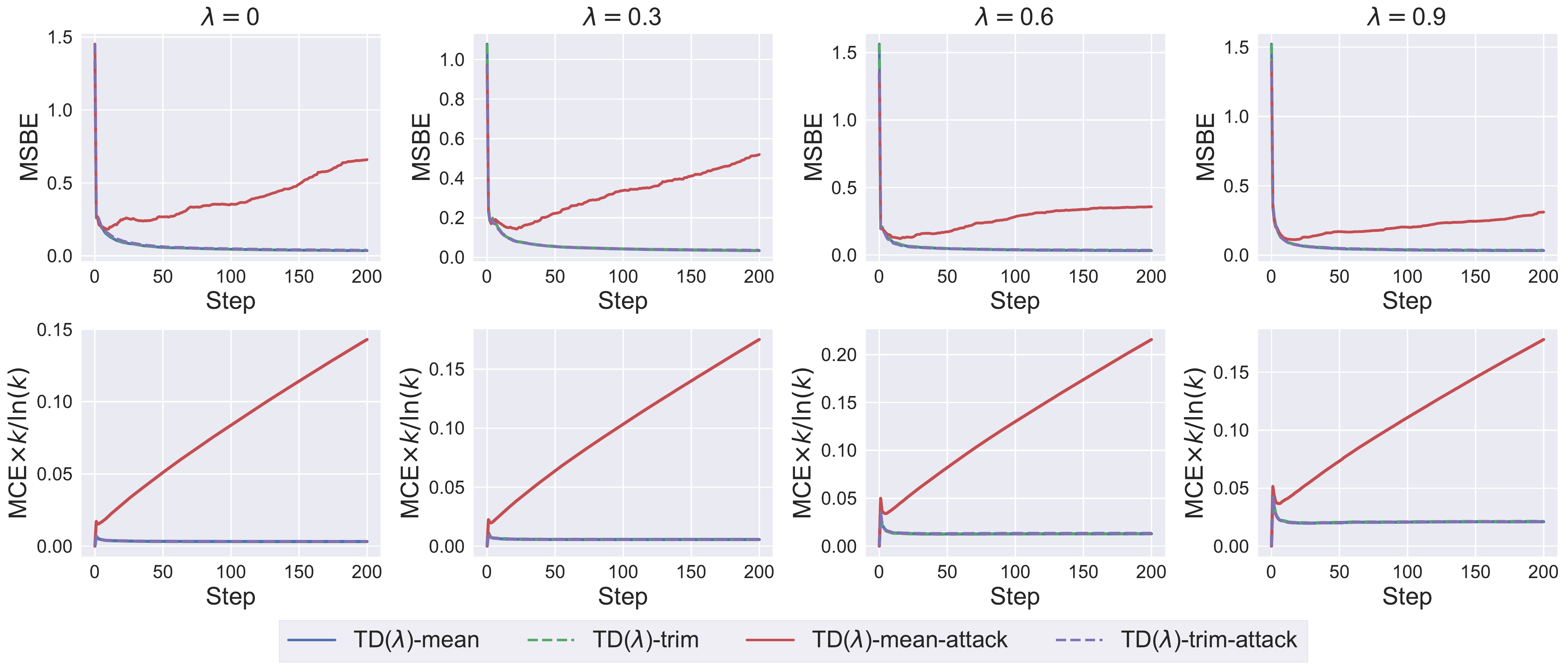}
	\vspace*{-0.3cm}
	\caption{\texttt{MSBE} and \texttt{MCE} under Gaussian noise attacks in Erdos-Renyi graph.}
	\label{fig:renyi-a1}
	\vspace*{-0.3cm}
\end{figure*}

 \textbf{Different Byzantine attacks.} The robustness of the algorithms has been tested under different attacks. In Figures \ref{fig:renyi-a0} and \ref{fig:renyi-a1}, we test same value and Gaussian noise attacks in the Erdos-Renyi graph. With \textit{same value} attacks, Byzantine agents always share zero vectors. In \textit{Gaussian noise} attacks, each Byzantine agent shares a randomly picked honest agent's parameter, polluted with Gaussian noise (with mean $0$ and standard deviation $1$). For $\lambda=0, 0.3, 0.6$, the step size is set as $\eta^k=\frac{0.1}{\sqrt{k}}$. For $\lambda=0.9$, $\eta^k=\frac{0.05}{\sqrt{k}}$. A total of 5 random Erdos-Renyi graphs are generated to calculate the averaged \texttt{MSBE} and \texttt{MCE}. In these experiments, $q_n$ is set to be the number of Byzantine agents in the graph.


\section{Conclusions}

\vspace{-0.8em}

This paper deals with temporal-difference learning in the MARL problem and proposes a novel \AlgorithmFullName. We establish the finite-time convergence rate and asymptotic learning error in the presence of Byzantine attacks. We show that the asymptotic learning error is determined by the variation of local rewards and the degree of network unsaturation. Numerical experiments corroborate the robustness of the proposed algorithm. Future work includes adaptively estimating the upper bound on the number of Byzantine agents.





%
%
%

\bibliography{abrv, ByzantineRL}

\newpage
\appendix
\begin{center}
	\Large \textbf{Supplementary Material} \\
\end{center}

\section{Convergence of Centralized TD($\lambda$)}
\label{sec:supportingLemma-convergence-rate-of-TD}

In these section, we will show the convergence rate of centralized TD($\lambda$).





To make the proof self-contained, we give the explicit expressions of $\mathbf{A}^*$ and $\bm{\bar{b}}^*$ defined in \eqref{eq:Ab}; see \cite{tsitsiklis1997analysis}. Let $\{s^{k}\}$ be the Markov chain defined by the underlying MRP, where $s^{k}\in\mathcal{S}$ has a transition probability matrix $\mathbf{P}\in\mathbb{R}^{|\mathcal{S}|\times|\mathcal{S}|}$, with $(s, s')$-th entry being $\mathcal{P}(s, s')$. We have
\begin{align}
	\mathbf{A}^* := \mathbf{\Phi}^{\top}\mathbf{D}(\mathbf{U}-\mathbf{I})\mathbf{\Phi}, \quad
	\bm{\bar{b}}^* := \mathbf{\Phi}^{\top}\mathbf{D}\sum_{\kappa=0}^{\infty}\lambda^\kappa(\gamma\mathbf{P})^{\kappa+1}\bm{r}^*,
\end{align}
where $\mathbf{\Phi}$, $\mathbf{D}$ and $\mathbf{U}$ are defined by
\begin{align}
	\mathbf{\Phi} &:= [\bm{\phi}(1),\bm{\phi}(2),\cdots,\bm{\phi}(s),\cdots,\bm{\phi}(|\mathcal{S}|)]^{\top}
	\in\mathbb{R}^{|\mathcal{S}|\times D}, \\
	\mathbf{D} &:= \text{diag}(\{\rho(s) | s\in\mathcal{S}\}), \\
	\mathbf{U} &:=(1-\lambda)\sum_{\kappa=0}^{\infty}\lambda^\kappa(\gamma\mathbf{P})^{\kappa+1}
	\in\mathbb{R}^{|\mathcal{S}|\times |\mathcal{S}|},
\end{align}
and the $s$-th element of vector $\bm{r}^*\in\mathbb{R}^{|\mathcal{S}|}$ is defined by
\begin{align}
	[\bm{r}^*]_s &:=
	\frac{1}{N}\sum_{n\in\mathcal{N}}\sum_{s'\in\mathcal{S}} \mathcal{P}(s, s')\mathcal{R}_n(s, s').
\end{align}



As we have discussed in Section \ref{sec:stationaryTD}, one key challenge of analyzing TD($\lambda$) is its biased update. To overcome this challenge, we turn to the concept of mixing time of Markov chain.
\begin{Definition}
	Given a constant $\eta>0$, we denote by $\tau(\eta)$ the mixing time of the Markov chain, given by
	\begin{align}
		\|\mathbb{E}[\mathbf{A}^k-\mathbf{A}^*|\mathcal{F}^{k-\tau(\eta)}]\| &\le \eta, ~ \forall k\ge\tau(\eta), \\
		\|\mathbb{E}[\bar{\bm{b}}^k-\bar{\bm{b}}^*|\mathcal{F}^{k-\tau(\eta)}]\| &\le \eta, ~ \forall k\ge\tau(\eta).
	\end{align}
\end{Definition}
Assumption \ref{assumption:markovChain} guarantees the Markov chain mixes at a geometric mixing rate, which means that there exists a constant $C$ such that for any given small constant $\eta$, $\tau(\eta)$ is bounded by
\begin{align}
	\tau(\eta) \le C\ln({1}/{\eta}).
\end{align}

In the following theorem,
we restate the convergence results derived in \cite{doan2019finite}.

\begin{Theorem}
	\label{theorem:convergence-rate-Central-TD-itretative}
	Suppose Assumptions \ref{assumption:markovChain} and \ref{assumption:NormalizedFeature} hold. If we use step size $\eta^k=\frac{\eta}{k+k^0}$ with large enough $k^0$, then given $\bm{\theta}$, the error after taking one step TD($\lambda$) follows from the two cases.
	
	Case 1. When $k\ge\tau(\eta^k)$, it holds
	\begin{align}
		\label{inequality:convergence-rate-Central-TD-itretative}
		\mathbb{E}[\|\bm{\theta}+\eta^k\bar g^k(\bm{\theta})-\bm{\theta}^{\infty}_{\lambda}\|^2]
		\le(1-\sigma_{\rm min}\eta^k)\|\bm{\theta}-\bm{\theta}^{\infty}_{\lambda}\|^2+2C_6 (\eta^k)^2 C \ln(\frac{k+1}{\eta}).
	\end{align}
	
	Case 2. When $k<\tau(\eta^k)$, it holds
	\begin{align}
		\|\bm{\theta}+\eta^k\bar g^k(\bm{\theta})-\bm{\theta}^{\infty}_{\lambda}\|^2
		\le  2\left(1+\eta^k\left(\frac{1+\gamma}{1-\gamma\lambda}\right)\right)^2
		\|\bm{\theta}-\bm{\theta}^{\infty}_{\lambda}\|^2
		+2(\eta^k)^2\left(\left(\frac{1+\gamma}{1-\gamma\lambda}\right)\| \bm{\theta}^{\infty}_{\lambda}\|+\left(\frac{R_{max}}{1-\gamma\lambda}\right)\right)^2,
	\end{align}
	where $R_{max}:=\frac{1}{N}\max_{s,s'\in\mathcal{S}}\sum_{n\in\mathcal{N}}\mathcal{R}(s, s')$ is the maximum reward, and the constants are defined as
	\begin{align}
		C_5 &:= 4\left(36+\frac{(229+42R_{max})(1+\gamma)^2}{(1-\gamma\lambda)^2}\right), \\
		C_6 &:= \|\bm{\theta}^{\infty}_{\lambda}\|^2C_5
		+2\left(
		32R_{max}^2+2\|\bm{\theta}^{\infty}_{\lambda}\|^2+1
		\right) \nonumber\\
		&~~~~~~+\frac{4\big(25R_{max}^2+16(R_{max}+1)^3+50(R_{max}+\|\bm{\theta}^{\infty}_{\lambda}\|^2)^2\big)(1+\gamma)^2}{(1-\gamma\lambda)^2}.
	\end{align}
\end{Theorem}

\section{TD($\lambda$) over a time-varying directed network}
\label{sec:time-varying-main-result-appendix}
We have showed the relation between Byrd-TD($\lambda$) and TD($\lambda$) over time-varying graph in section \ref{sec:time-varying-main-result}. Before we prove Theorem \ref{theorem:consensu-rate-theorem-global-time-varying} and \ref{theorem:convergence-trimmedMean-time-varying}, we define some auxiliary variables.

\subsection{Construction of auxiliary variable $\hat{\mathbf{\Theta}}^{k}$}

The lack of consensus of $\bm{\theta}^{k+1}_n$ in the initial stage brings difficulty to our proof. Therefore, we look for a variable $\bm{\omega}^k\in\mathbb{R}^D$, for which all $\bm{\theta}^{k+1}_n$ will converge to.

To gain insights, we expand \eqref{equation:update-trimmedMean-matrix} as
\begin{align}
	\label{equation:update-trimmedMean-matrix-fromzero}
	\mathbf{\Theta}^{k}_d
	=&\mathbf{Y}^{k-1}(d)\mathbf{\Theta}^{k-1}_d+\eta^{k-1} \mathbf{G}^{k-1}_d(\mathbf{\Theta}^k) \nonumber\\
	=&\mathbf{Y}^{k-1}(d)\mathbf{Y}^{k-2}(d)\cdots\mathbf{Y}^0(d)\mathbf{\Theta}^{0}_d
	+\sum_{\kappa=0}^{k-1}\mathbf{Y}^{k-1}(d)\mathbf{Y}^{k-2}(d)\cdots\mathbf{Y}^{\kappa+1}(d)\eta^\kappa \mathbf{G}^\kappa_d(\mathbf{\Theta}^\kappa) \nonumber\\
	=&\bar{\mathbf{Y}}_0^{k-1}(d)\mathbf{\Theta}^{0}_d
	+\sum_{\kappa=0}^{k-1}\bar{\mathbf{Y}}_{\kappa+1}^{k-1}(d)\eta^\kappa \mathbf{G}^\kappa_d(\mathbf{\Theta}^\kappa),
\end{align}
where the transition matrix $\bar{\mathbf{Y}}_{k_0}^k(d)$ is defined as the product of a series of $\mathbf{Y}^k(d)$, as
\begin{align}
	\label{definition:Psi-k}
	\bar{\mathbf{Y}}_{k_0}^k(d):=\mathbf{Y}^k(d)\mathbf{Y}^{k-1}(d)\cdots \mathbf{Y}^{k_0}(d).
\end{align}
The following lemma shows the limit of $\bar{\mathbf{Y}}_{k_0}^k(d)$ exists.
\begin{Lemma}
	\label{lemma:limit-of-product-of-Yk}
	For any possible series $\{\mathbf{Y}^k(d)\}$ satisfying conditions (c1)--(c6), the limit of $\bar{\mathbf{Y}}_{k_0}^k(d)$ exists, denoted as $\bm{1}\bm{p}_{k_0}^{\top}(d)$, and $\bar{\mathbf{Y}}_{k_0}^k(d)$ converges to its limit at linear rate, such that
	\begin{align}
		\lim_{k\to\infty}\bar{\mathbf{Y}}_{k_0}^k(d) &= \bm{1}\bm{p}_{k_0}^{\top}(d) 
		\label{equation:lemma-Y-bar-limit} \\
		\|\bar{\mathbf{Y}}_{k_0}^k(d) - \bm{1}\bm{p}_{k_0}^{\top}(d)\|^2 
		&\le 4N^2\mu^{k-k_0+1}
		\label{equation:lemma-Y-bar-limit-rate} 
	\end{align}
	where $\bm{p}_{k_0}(d)$ is a stochastic vector and all elements in $\bm{p}_{k_0}(d)$ are within $\left[0, \frac{1}{\min_{n\in\mathcal{N}}\{N_n+B_n-2q_n+1\}}\right]$.
\end{Lemma}


With Lemma \ref{lemma:limit-of-product-of-Yk}, we introduce a matrix $\hat{\mathbf{\Theta}}^{k} := \bm{1}(\bm{\omega}^{k})^{\top}$ and a vector $\hat{\mathbf{\Theta}}^{k}_d := [\bm{\omega}^{k}]_d\bm{1}$, where $[\bm{\omega}^{k}]_d$ is the $d$-th element of $\bm{\omega}^k\in\mathbb{R}^D$,
The vector $\hat{\mathbf{\Theta}}^{k}_d$ asymptotically approximates the sequence $\mathbf{\Theta}^{k}_d$, and satisfies the following recursion
\begin{align}
	\label{definition:consensu-mapping}
	[\bm{\omega}^{k}]_d :=& \bm{p}_0^{\top}(d)\hat{\mathbf{\Theta}}^{0}_d+\sum_{\kappa=0}^{k-1}\bm{p}^{\top}_{\kappa+1}(d)\eta^\kappa \mathbf{G}^\kappa_d(\hat{\mathbf{\Theta}}^{\kappa}).
\end{align}

With the auxiliary variable $\bm{\omega}^k$ and $\hat{\mathbf{\Theta}}^{k}$, we can begin our proofs of Theorems \ref{theorem:consensu-rate-theorem-global-time-varying} and \ref{theorem:convergence-trimmedMean-time-varying}.

\subsection{Proof of Theorem \ref{theorem:consensu-rate-theorem-global-time-varying}}
Because $\|\bm{\theta}^{k}_n-\bar{\bm{\theta}}^{k}\|^2$ can be bounded by $\|\mathbf{\Theta}^{k}-\hat{\mathbf{\Theta}}^{k}\|^2_F$ as
\begin{align}
	\frac{1}{N}\sum_{n\in\mathcal{N}}\|\bm{\theta}^{k}_n-\bar{\bm{\theta}}^{k}\|^2
	=&	\frac{1}{N}\sum_{n\in\mathcal{N}}\|\bm{\theta}^{k}_n-\bm{\omega}^k\|^2
	- \|\bm{\omega}^k-\bar{\bm{\theta}}^{k}\|^2
	\le \frac{1}{N} \|\mathbf{\Theta}^{k}-\hat{\mathbf{\Theta}}^{k}\|^2_F,
\end{align}
we will bound $\|\bm{\theta}^{k}_n-\bar{\bm{\theta}}^{k}\|^2$ by analyzing $\|\mathbf{\Theta}^{k}-\hat{\mathbf{\Theta}}^{k}\|^2_F$. We first give the bound of difference as
\begin{align}
	\label{inequality:different-to-consensus-mapping}
	&\|\mathbf{\Theta}^{k}-\hat{\mathbf{\Theta}}^{k}\|^2_F
	= \sum_{d=1}^{D}\|\mathbf{\Theta}^{k}_d-\hat{\mathbf{\Theta}}^{k}_d\|^2
	\nonumber \\
	\le& 2\sum_{d=1}^{D}\|
	(\bar{\mathbf{Y}}_0^{k-1}(d)-\bm{1}\bm{p}_0^{\top}(d))\mathbf{\Theta}^{0}_d\|^2
	+2\sum_{d=1}^{D}\sum_{\kappa=0}^{k-1}\|\left(\bar{\mathbf{Y}}_{\kappa+1}^{k-1}(d)-\bm{1}\bm{p}_{\kappa+1}^{\top}(d)\right)\eta^\kappa \mathbf{G}^\kappa_d(\mathbf{\Theta}^\kappa)\|^2,
\end{align}
which comes from inequality $\|x+y\|^2\le2\|x\|^2+2\|y\|^2$.

Since both $\bar{\mathbf{Y}}_{\kappa}^{k-1}(d)$ and $\bm{1}\bm{p}_{\kappa}^{\top}(d)$ are row-stochastic matrices, we have $\left(\bar{\mathbf{Y}}_{\kappa}^{k-1}(d)-\bm{1}\bm{p}_{\kappa}^{\top}(d)\right)\bm{1}=\bm{0}$, which implies that the first term of \eqref{inequality:different-to-consensus-mapping} can be bounded by
\begin{align}
	\label{inequality:different-to-consensus-mapping-1}
	&2\sum_{d=1}^{D}\|
	(\bar{\mathbf{Y}}_0^{k-1}(d)-\bm{1}\bm{p}_0^{\top}(d))\mathbf{\Theta}^{0}_d\|^2
	\nonumber \\
	=& 2\sum_{d=1}^{D}\|
	(\bar{\mathbf{Y}}_0^{k-1}(d)-\bm{1}\bm{p}_0^{\top}(d))(\mathbf{\Theta}^{0}_d-\bm{\omega}^0_d\bm{1})\|^2
	\nonumber \\
	\le& 2\sum_{d=1}^{D}
	\|\bar{\mathbf{Y}}_0^{k-1}(d)-\bm{1}\bm{p}_0^{\top}(d)\|^2
	\|\mathbf{\Theta}^{0}_d-\bm{\omega}^0_d\bm{1}\|^2
	\nonumber \\
	\le& 8\sum_{d=1}^{D} N^2\mu^{k}
	\|\mathbf{\Theta}^{0}_d-\bm{\omega}^0_d\bm{1}\|^2
	= 0,
\end{align}
where the last equality is from $\hat{\mathbf{\Theta}}^{0} = \mathbf{\Theta}^{0}$.

Noticing the inequality that
\begin{align}
	&\|\left(\bar{\mathbf{Y}}_{\kappa+1}^{k-1}(d)-\bm{1}\bm{p}_{\kappa+1}^{\top}(d)\right)\eta^\kappa \mathbf{G}^\kappa_d(\mathbf{\Theta}^\kappa)\|^2
	\nonumber \\
	=& \|\left(\bar{\mathbf{Y}}_{\kappa+1}^{k-1}(d)-\bm{1}\bm{p}_{\kappa+1}^{\top}(d)\right)\eta^\kappa (\mathbf{G}^\kappa_d(\mathbf{\Theta}^\kappa)-[\bm{\bar {g}}^\kappa(\bm{\omega}^\kappa)]_d\bm{1})\|^2
	\nonumber \\
	\le& \|\bar{\mathbf{Y}}_{\kappa+1}^{k-1}(d)-\bm{1}\bm{p}_{\kappa+1}^{\top}(d)\|^2
	\|\eta^\kappa (\mathbf{G}^\kappa_d(\mathbf{\Theta}^\kappa)-[\bm{\bar {g}}^\kappa(\bm{\omega}^\kappa)]_d\bm{1})\|^2
	\nonumber \\
	\le& 4N^2\mu^{k-\kappa-1}\|\eta^\kappa (\mathbf{G}^\kappa_d(\mathbf{\Theta}^\kappa)-[\bm{\bar {g}}^\kappa(\bm{\omega}^\kappa)]_d\bm{1})\|^2,
\end{align}
the second term of \eqref{inequality:different-to-consensus-mapping} can be bounded by
\begin{align}
	\label{inequality:different-to-consensus-mapping-2}
	&2\sum_{d=1}^{D}\sum_{\kappa=0}^{k-1}\|\left(\bar{\mathbf{Y}}_{\kappa+1}^{k-1}(d)-\bm{1}\bm{p}_{\kappa+1}^{\top}(d)\right)\eta^\kappa \mathbf{G}^\kappa_d(\mathbf{\Theta}^\kappa)\|^2
	\nonumber \\
	\le& 8\sum_{\kappa=0}^{k-1}N^2
	\mu^{k-\kappa-1}
	\sum_{d=1}^{D}\|\eta^\kappa (\mathbf{G}^\kappa_d(\mathbf{\Theta}^\kappa)-[\bm{\bar {g}}^\kappa(\bm{\omega}^\kappa)]_d\bm{1})\|^2
	\nonumber \\
	=& 8N^2\sum_{\kappa=0}^{k-1}
	\mu^{k-\kappa-1}
	(\eta^\kappa)^2\|\mathbf{G}^\kappa(\mathbf{\Theta}^\kappa)-\bm{1}\bm{\bar {g}}^\kappa(\bm{\omega}^\kappa)^{\top}\|^2_F.
\end{align}

Note that the right-hand side (RHS) of \eqref{inequality:different-to-consensus-mapping-2} is the difference between the local and global increments, which will be bounded in the following lemma.

\begin{Lemma}
	\label{lemma:different-local-global-update}
	For any $\bm{\theta}_n$ and $\bm{\theta}'$, the difference between local and global increment is upper bounded by
	\begin{align}
		\label{inequatliy:different-local-global-update}
		\frac{1}{N}\sum_{n\in\mathcal{N}}\|\bm{g}^k_n(\bm{\theta}_n)-\bm{\bar {g}}^k(\bm{\theta}')\|^2
		\le 2\left(\frac{1+\gamma}{1-\gamma\lambda}\right)^2\frac{1}{N}\sum_{n\in\mathcal{N}}\|\bm{\theta}_n-\bm{\theta}'\|^2
		+ \frac{2\delta^2}{(1-\gamma\lambda)^2}.
	\end{align}
\end{Lemma}

The proof of Lemma \ref{lemma:different-local-global-update} is left to the supplementary, as well as the full version of this paper \cite{wu2020byzantine}. With Lemma \ref{lemma:different-local-global-update},
the last term of \eqref{inequality:different-to-consensus-mapping-2} can be bounded further
\begin{align}
	\label{inequality:different-to-consensus-mapping-2-2}
	&\|\mathbf{G}^\kappa(\mathbf{\Theta}^\kappa)-\bm{1}\bm{\bar {g}}^\kappa(\bm{\omega}^\kappa)^{\top}\|^2_F
	=\sum_{n\in\mathcal{N}}\|\bm{g}^\kappa_n(\bm{\theta}^\kappa_n)-\bm{\bar {g}}^\kappa(\bm{\omega}^\kappa)\|^2\nonumber\\
	\le& 2\left(\frac{1+\gamma}{1-\gamma\lambda}\right)^2\sum_{n\in\mathcal{N}}\|\bm{\theta}^\kappa_n-\bm{\omega}^\kappa\|^2
	+ \frac{2N\delta^2}{(1-\gamma\lambda)^2}
	\nonumber \\
	=& 2\left(\frac{1+\gamma}{1-\gamma\lambda}\right)^2\|\mathbf{\Theta}^{\kappa}-\hat{\mathbf{\Theta}}^{\kappa}\|^2_F
	+ \frac{2N\delta^2}{(1-\gamma\lambda)^2}.
\end{align}

Plugging \eqref{inequality:different-to-consensus-mapping-1}, \eqref{inequality:different-to-consensus-mapping-2} and \eqref{inequality:different-to-consensus-mapping-2-2} into \eqref{inequality:different-to-consensus-mapping}, we have
\begin{align}
	\label{inequality:different-to-consensus-mapping-iterative-1}
	\|&\mathbf{\Theta}^{k}-\hat{\mathbf{\Theta}}^{k}\|^2_F 
	=\frac{16N^3\delta^2}{(1-\gamma\lambda)^2}
	\sum_{\kappa=0}^{k-1}\mu^{k-\kappa-1}(\eta^\kappa)^2
	+16N^2(\frac{1+\gamma}{1-\gamma\lambda})^2\sum_{\kappa=0}^{k-1}
	\mu^{k-\kappa-1}
	(\eta^\kappa)^2\|\mathbf{\Theta}^{\kappa}-\hat{\mathbf{\Theta}}^{\kappa}\|^2_F.
\end{align}

We can bound the two terms in the RHS of \eqref{inequality:different-to-consensus-mapping-iterative-1} by using the following lemma.
\begin{Lemma}
	\label{lemma:bound-consensus-auxiliary}
	When the step size $\eta^k$ in Algorithm \ref{algorithm:TD_TrimmedMean} satisfies 
	\begin{align}
		1\le
		\left(\frac{\eta^{k-1}}{\eta^{k}}\right)^2
		\le \frac{2}{1+\alpha_1},
	\end{align}
	and for some $\alpha_1\in(0, 1)$ and $\alpha_2\ge 0$, the iterates $\{x^k\}$ satisfy
	\begin{align}
		x^{k+1}\le \alpha_1 x^k + \alpha_2(\eta^k)^2
		~~\text{and}~~
		x^0 \le \alpha_2(\eta^0)^2    ,
	\end{align}
	then $x^k$ has an upper bound
	\begin{align}
		x^k\le \frac{2\alpha_2}{1-\alpha_1}(\eta^k)^2.
	\end{align}
\end{Lemma}
Lemma \ref{lemma:bound-consensus-auxiliary} can be proved by induction.

In order to bound the first term in the RHS of \eqref{inequality:different-to-consensus-mapping-iterative-1}, we define an auxiliary variable $x_1^k$ as
\begin{align}
	x_1^k:= \sum_{\kappa=0}^{k-1}\mu^{k-\kappa-1}(\eta^\kappa)^2,
\end{align}
which has relation $x_1^{k+1} = \mu x_1^k+(\eta^k)^2$. Therefore, Lemma \ref{lemma:bound-consensus-auxiliary} implies that
\begin{align}
	\label{inequality:different-to-consensus-mapping-iterative-1-1}
	\sum_{\kappa=0}^{k-1}\mu^{k-\kappa-1}(\eta^\kappa)^2
	= x^k_1 \le \frac{2}{1-\mu}(\eta^k)^2.
\end{align}
Substituting \eqref{inequality:different-to-consensus-mapping-iterative-1-1} into \eqref{inequality:different-to-consensus-mapping-iterative-1} yields that
\begin{align}
	\label{inequality:different-to-consensus-mapping-iterative-2}
	\|&\mathbf{\Theta}^{k}-\hat{\mathbf{\Theta}}^{k}\|^2_F
	\le \frac{32N^3\delta^2}{(1-\gamma\lambda)^2(1-\mu)} (\eta^k)^2
	+16N^2(\frac{1+\gamma}{1-\gamma\lambda})^2\sum_{\kappa=0}^{k-1}
	\mu^{k-\kappa-1}
	(\eta^\kappa)^2\|\mathbf{\Theta}^{\kappa}-\hat{\mathbf{\Theta}}^{\kappa}\|^2_F.
\end{align}
Similarly, to bound the second terms in the RHS of \eqref{inequality:different-to-consensus-mapping-iterative-2}, we define another auxiliary variable $x^k_2$ that
\begin{align}
	x^k_2& := \frac{32N^3\delta^2}{(1-\gamma\lambda)^2(1-\mu)}(\eta^k)^2
	+16N^2(\frac{1+\gamma}{1-\gamma\lambda})^2\sum_{\kappa=0}^{k-1}
	\mu^{k-\kappa-1}
	(\eta^\kappa)^2\|\mathbf{\Theta}^{\kappa}-\hat{\mathbf{\Theta}}^{\kappa}\|^2_F,
\end{align}
then we have $\|\mathbf{\Theta}^{k}-\hat{\mathbf{\Theta}}^{k}\|^2_F\le x^k_2$ and
\begin{align}
	x^{k+1}_2 \le& \frac{32N^3\delta^2}{(1-\gamma\lambda)^2(1-\mu)}(\eta^{k+1})^2
	+16N^2(\frac{1+\gamma}{1-\gamma\lambda})^2\sum_{\kappa=0}^{k}
	\mu^{k-\kappa-1}
	(\eta^\kappa)^2 x^\kappa_2
	\nonumber \\
	=& \left(\mu + 16N^2(\frac{1+\gamma}{1-\gamma\lambda})^2 (\eta^k)^2\right) x^k_2
	+ \frac{32N^3\delta^2}{(1-\gamma\lambda)^2(1-\mu)}
	\left((\eta^{k+1})^2-\mu(\eta^{k})^2\right)
	\nonumber \\
	\le& \frac{1+\mu}{2} x^k_2
	+ \frac{32N^3\delta^2}{(1-\gamma\lambda)^2}(\eta^{k})^2,
\end{align}
where the last inequality comes from the step size condition \eqref{condition:stepsize-consensus-1}. Applying Lemma \ref{lemma:bound-consensus-auxiliary} again yields that
\begin{align}
	\label{inequality:bound-Theta-consensus}
	\|\mathbf{\Theta}^{k}-\hat{\mathbf{\Theta}}^{k}\|^2_F
	\le x^k_2
	\le \frac{128N^3\delta^2}{(1-\gamma\lambda)^2(1-\mu)}(\eta^{k})^2.
\end{align}

Now we reach the following inequality
\begin{align}
	\frac{1}{N}\sum_{n\in\mathcal{N}}\|\bm{\theta}^{k}_n-\bar{\bm{\theta}}^{k}\|^2
	\!\le\! \frac{1}{N} \|\mathbf{\Theta}^{k}-\hat{\mathbf{\Theta}}^{k}\|^2_F
	\le\frac{C_1\mu_{\mathcal{G}}(\eta^{k})^2}{2},
\end{align}
which completes the proof.

\subsection{Proof of Theorem \ref{theorem:convergence-trimmedMean-time-varying}}
By analyzing $\|\bm{\theta}^{k+1}_n-\bm{\theta^{\infty}}_{\lambda}\|^2$, we find that
\begin{align}\label{pf.thm5_1}
	& \frac{1}{N}\sum_{n\in\mathcal{N}}\mathbb{E}\|\bm{\theta}^{k+1}_n-\bm{\theta^{\infty}}_{\lambda}\|^2
	\nonumber\\
	\le& \frac{1}{N}\sum_{n\in\mathcal{N}}\left(2\mathbb{E}\|\bm{\theta}^{k+1}_n-\bm{\omega}^{k+1}\|^2+2\mathbb{E}\|\bm{\omega}^{k+1}-\bm{\theta^{\infty}}_{\lambda}\|^2\right)
	\nonumber\\
	=& \frac{2}{N}\mathbb{E}\|\mathbf{\Theta}^{k+1}-\hat{\mathbf{\Theta}}^{k+1}\|^2_F+2\mathbb{E}\|\bm{\omega}^{k+1}-\bm{\theta^{\infty}}_{\lambda}\|^2.
\end{align}

Since $\|\mathbf{\Theta}^{k+1}-\hat{\mathbf{\Theta}}^{k+1}\|^2_F$ has been bounded by \eqref{inequality:bound-Theta-consensus}, we will proceed with $\|\bm{\omega}^{k+1}-\bm{\theta^{\infty}}_{\lambda}\|^2$. To do so, note that $[\bm{\omega}^{k}]_d$ in \eqref{definition:consensu-mapping} can be expressed as
\begin{align}
	\label{equation:update-consensus-mapping}
	[\bm{\omega}^{k+1}]_d =& [\bm{\omega}^{k}]_d+\bm{p}_{k+1}^{\top}(d)\eta^k \mathbf{G}^k_d(\hat{\mathbf{\Theta}}^{k}),
\end{align}

Observe that the update of $[\bm{\omega}^{k}]_d$ has a similar form of the TD($\lambda$) update, so we manage to analyze its convergence using the techniques in the TD($\lambda$) analysis. However, TD($\lambda$) over the decentralized and time-varying graph cannot reach a minimum of $F(\bm{\theta})$ like the original TD($\lambda$) since $\bm{p}_{k+1}(d) \neq \frac{1}{N}\bm{1}$ usually. This is our main challenge, comparing to the analysis in \cite{doan2019finite}.


Stacking \eqref{equation:update-consensus-mapping} over all dimensions $d$, we have
\begin{align}
	&\|\bm{\omega}^{k+1}-\bm{\theta^{\infty}}_{\lambda}\|^2 
	= \sum_{d=1}^{D}\Big([\bm{\omega}^{k+1}]_d-[\bm{\theta^{\infty}}_{\lambda}]_d\Big)^2
	\sum_{d=1}^{D}\Big([\bm{\omega}^{k}]_d+\bm{p}_{k+1}^{\top}(d)\eta^k \mathbf{G}^k_d(\hat{\mathbf{\Theta}}^{k})-[\bm{\theta^{\infty}}_{\lambda}]_d\Big)^2,
\end{align}
and using the following inequality
\begin{align}
	&\Big([\bm{\omega}^{k}]_d+\bm{p}_{k+1}^{\top}(d)\eta^k \mathbf{G}^k_d(\hat{\mathbf{\Theta}}^{k})-[\bm{\theta^{\infty}}_{\lambda}]_d\Big)^2
	\nonumber \\
	\le&\frac{1}{1-\beta^k}\Big([\bm{\omega}^{k}]_d
	+\frac{1}{N}\bm{1}^{\top}\eta^k \mathbf{G}^k_d(\hat{\mathbf{\Theta}}^{k})
	-[\bm{\theta^{\infty}}_{\lambda}]_d\Big)^2
	+\frac{1}{\beta^k}\Big((\bm{p}_{k+1}(d)-\frac{1}{N}\bm{1})^{\top}\eta^k \mathbf{G}^k_d(\hat{\mathbf{\Theta}}^{k})\Big)^2,
\end{align}
for any $0<\beta^k<1$, we have
\begin{align}
	\label{inequality:convergence-rate-consensus-mapping-1}
	\|\bm{\omega}^{k+1}-\bm{\theta^{\infty}}_{\lambda}\|^2
	&	= \frac{1}{1-\beta^k}\|\bm{\omega}^{k}
	+\eta^k\bm{\bar g}^k(\bm{\omega}^{k})
	-\bm{\theta^{\infty}}_{\lambda}\|^2
	+\frac{1}{\beta^k}\sum_{d=1}^{D}\Big((\bm{p}_{k+1}(d)-\frac{1}{N}\bm{1})^{\top}\eta^k \mathbf{G}^k_d(\hat{\mathbf{\Theta}}^{k})\Big)^2.
\end{align}



We now bound the second term at the RHS of \eqref{inequality:convergence-rate-consensus-mapping-1}. With Cauchy-Schwarz inequality, it holds
\begin{align}
	\label{inequality:convergence-rate-consensus-mapping-1-1}
	&\sum_{d=1}^{D}\Big((\bm{p}_{k+1}(d)-\frac{1}{N}\bm{1})^{\top}\eta^k \mathbf{G}^k_d(\hat{\mathbf{\Theta}}^{k})\Big)^2
	\\
	=&\sum_{d=1}^{D}
	\Big((\bm{p}_{k+1}(d)-\frac{1}{N}\bm{1})^{\top}\eta^k
	\left(\mathbf{G}^k_d(\hat{\mathbf{\Theta}}^{k})-[\bm{\bar {g}}^k(\bm{\omega}^k)]_d\bm{1}\right)
	\Big)^2
	\nonumber \\
	\le & (\eta^k)^2\sum_{d=1}^{D}\|\bm{p}_{k+1}(d)-\frac{1}{N}\bm{1}\|^2
	\|\mathbf{G}^k_d(\hat{\mathbf{\Theta}}^{k})-[\bm{\bar {g}}^k(\bm{\omega}^k)]_d\bm{1}\|^2,\nonumber 
\end{align}
where the equality comes from the fact that 
\begin{equation}
	(\bm{p}_{k+1}(d)-\frac{1}{N}\bm{1})^{\top}\bm{1}=0.
\end{equation}
From Lemma \ref{lemma:limit-of-product-of-Yk}, it follows that $\bm{p}_{k_0}(d)$ is a stochastic vector and the elements of $\bm{p}_{k_0}(d)$ are 
\begin{equation}
	\bm{p}_{k_0}(d)\in \left[0, \frac{\bm{1}}{\min_{n\in\mathcal{N}}\{N_n+B_n-2q_n+1\}}\right].
\end{equation}

Therefore, the first term in \eqref{inequality:convergence-rate-consensus-mapping-1-1} has a uniform upper bound
\begin{equation}
	\Big\|\bm{p}_{k+1}(d)-\frac{\bm{1}}{N}\Big\|^2
	\!\!\le\! \frac{1}{\min_{n\in\mathcal{N}}\{N_n+B_n-2q_n+1\}}-\frac{1}{N}
	\!=\! \frac{D_{\mathcal{G}}}{N}.
\end{equation}
With this inequality, the RHS of \eqref{inequality:convergence-rate-consensus-mapping-1-1} can be bounded by
\begin{align}
	\label{inequality:convergence-rate-consensus-mapping-1-2}
	& (\eta^k)^2\sum_{d=1}^{D}\|\bm{p}_{k+1}(d)-\frac{1}{N}\bm{1}\|^2
	\left\|\mathbf{G}^k_d(\hat{\mathbf{\Theta}}^{k})-\bm{\bar {g}}_d^k(\bm{\omega}^k)\bm{1}\right\|^2
	\nonumber \\
	\le & (\eta^k)^2
	\frac{D_{\mathcal{G}}}{N}
	\sum_{d=1}^{D}
	\left\|\mathbf{G}^k_d(\hat{\mathbf{\Theta}}^{k})-\bm{\bar {g}}_d^k(\bm{\omega}^k)\bm{1}\right\|^2
	\nonumber \\
	= & (\eta^k)^2
	\frac{D_{\mathcal{G}}}{N} \|\mathbf{G}^k(\hat{\mathbf{\Theta}}^k)-\bm{1}\bm{\bar {g}}^k(\bm{\omega}^k)^{\top}\|^2_F
	\nonumber \\
	\le & (\eta^k)^2\frac{2D_{\mathcal{G}}\delta^2}{(1-\gamma\lambda)^2},
\end{align}
where $D_{\mathcal{G}}$ was defined in \eqref{definition:networkunsaturability} and the last inequality comes from Lemma \ref{lemma:different-local-global-update} that
\begin{align}
	\|\mathbf{G}^k(\hat{\mathbf{\Theta}}^k)-\bm{1}\bm{\bar {g}}^k(\bm{\omega}^k)^{\top}\|^2_F
	=&\sum_{n\in\mathcal{N}}\|\bm{ g}^k_n(\bm{\omega}^k)-\bm{\bar {g}}^k(\bm{\omega}^k)\|^2
	\le \frac{2N\delta^2}{(1-\gamma\lambda)^2}.
\end{align}
With \eqref{inequality:convergence-rate-consensus-mapping-1-1} and  \eqref{inequality:convergence-rate-consensus-mapping-1-2}, \eqref{inequality:convergence-rate-consensus-mapping-1} becomes
\begin{align}
	\label{inequality:convergence-rate-consensus-mapping-1.5}
	\|\bm{\omega}^{k+1}-\bm{\theta^{\infty}}_{\lambda}\|^2
	\le\frac{1}{1-\beta^k}\|\bm{\omega}^{k}
	+\eta^k\bm{\bar g}^k(\bm{\omega}^{k})
	-\bm{\theta^{\infty}}_{\lambda}\|^2
	+\frac{(\eta^k)^2}{\beta^k} \frac{2D_{\mathcal{G}}\delta^2}{(1-\gamma\lambda)^2}. \nonumber
\end{align}

The first term at the RHS of \eqref{inequality:convergence-rate-consensus-mapping-1.5} can be bounded in terms of $\mathbb{E}\|\bm{\omega}^{k}-\bm{\theta^{\infty}}_{\lambda}\|^2$ by Theorem \ref{theorem:convergence-rate-Central-TD-itretative}, and the following lemma bounds the sequence $\{\mathbb{E}\|\bm{\omega}^{k}-\bm{\theta^{\infty}}_{\lambda}\|^2\}$.

\begin{Lemma}
	\label{lemma:convergence-omega}
	If $\{\bm{\omega}^{k}\}$ is the sequence generated by \eqref{equation:update-consensus-mapping}, the expected optimal gap can be bounded by
	\begin{align}
		\mathbb{E}\|\bm{\omega}^{k}-\bm{\theta^{\infty}}_{\lambda}\|^2
		\le&\frac{C_2}{2}\left(\frac{1}{k+k^0}\right)^\epsilon
		\!\!+\frac{C_3}{2}\varphi_{\epsilon}(k)
		+\frac{C_4}{2} \frac{D_{\mathcal{G}}\delta^2}{(1-\gamma\lambda)^2}.
	\end{align}
\end{Lemma}

The proof of Lemma \ref{lemma:convergence-omega} is left to the supplementary, as well as the full version of this paper \cite{wu2020byzantine}. Together with Theorem \ref{theorem:consensu-rate-theorem-global-time-varying}, we can use \eqref{pf.thm5_1} to present the convergence rate of $\bm{\theta}^k_n$ as
\begin{align}
	&\frac{1}{N}\sum_{n\in\mathcal{N}}\mathbb{E}\|\bm{\theta}^{k+1}_n-\bm{\theta^{\infty}}_{\lambda}\|^2
	\nonumber\\
	\le& \frac{2}{N}\mathbb{E}\|\mathbf{\Theta}^{k+1}-\hat{\mathbf{\Theta}}^{k+1}\|^2_F+2\mathbb{E}\|\bm{\omega}^{k+1}-\bm{\theta^{\infty}}_{\lambda}\|^2
	\nonumber\\
	\le& C_1\mu_{\mathcal{G}}(\eta^{k})^2
	+\frac{C_2}{(k+k^0)^\epsilon}
	+ C_3\varphi_\epsilon(k)
	+ C_4\frac{D_{\mathcal{G}}\delta^2}{(1-\gamma\lambda)^2}.
\end{align}
Now we complete the proof.

\section{Consensus and Convergence of Byrd-TD($\lambda$)}
\label{sec:Consensus-Convergence-Byrd-TD}

We have shown that when the topology $\mathbf{Y}^k(d)$ satisfies Conditions (c1)--(c6), TD($\lambda$) over a time-varying directed networks guarantees both consensus and convergence.

Now we will show that Byrd-TD($\lambda$) can be described by \eqref{equation:update-trimmedMean-matrix} with $\mathbf{Y}^k(d)$ satisfying Conditions (c1)--(c6).
Therefore, Theorems \ref{theorem:consensu-rate-theorem-global-time-varying} and \ref{theorem:convergence-trimmedMean-time-varying} can be applied to Byrd-TD($\lambda$).

\subsection{Proofs of Theorems \ref{theorem:consensu-rate-theorem-global} and \ref{theorem:convergence-trimmedMean}}
Now the only thing we need to do for the proofs of Theorems \ref{theorem:consensu-rate-theorem-global} and \ref{theorem:convergence-trimmedMean} is to verify that our constructed $\mathbf{Y}^k(d)$ satisfies Conditions (c1)--(c6) under Assumption \ref{assumption:networkStructure}.

Conditions (c1), (c2) and (c3) directly follow the similar results in \cite[Claim 2]{vaidya2012matrix}. In addition, the proofs of (c4) and (c5) for Case 1 are also simple.

We only show (c4) and (c5) for Case 2 and (c6) here.

\noindent\textbf{Verification of (c4) for Case 2}. 

In Case 2, because $0\le 3q_n<N_n$, we have $B_n \le q_n < N_n+B_n-2q_n$ and then $B^{k*}_n(d)\le \min\{N_n+B_n-2q_n, B_n\} = B_n$. As a result, it follows
\begin{align}
	\frac{q_n-B_n}{N_n+B_n-2q_n-B^{k*}_n(d)}
	\le \frac{q_n-B_n}{N_n-2q_n}.
\end{align}

For $m\in \mathcal{N}_n^{k*}(d)\cap\mathcal{N}_n$, we have
\begin{align}
	[\mathbf{Y}^k(d)]_{nm} =& \left(
	1-\frac{q_n-B_n}{N_n+B_n-2q_n-B^{k*}_n(d)}
	\right)\frac{1}{N^*_n}
	\ge \frac{N_n+B_n-3q_n}{N_n-2q_n}\frac{1}{N^*_n}.
\end{align}



For $m\in\mathcal{L}_n^{k+}(d)$ (or $\mathcal{L}_n^{k-}(d)$), there exists $m'\in\mathcal{L}_n^{k-}(d)$ (or $\mathcal{L}_n^{k+}(d)$), such that $y(n',m)+y(n',m')=1$ and therefore at least one among $y(n',m)$ and $y(n',m')$ is no less than $\frac{1}{2}$. Then at least one among $\sum_{n'\in\mathcal{N}_n^{k*}(d)}y(n',m)$ and $\sum_{n'\in\mathcal{N}_n^{k*}(d)}y(n',m')$ is no less than $|\mathcal{N}_n^{k*}(d)|/4$. 

Without loss of generality, assuming $\sum_{n'\in\mathcal{N}_n^{k*}(d)}y(n',m)$ is no less than $|\mathcal{N}_n^{k*}(d)|/4$, it holds
\begin{align}
	[\mathbf{Y}^k(d)]_{nm}
	=&\frac{1}{q_n-B_n+B^{k*}_n(d)}\sum_{n'\in\mathcal{N}_n^{k*}(d)\cap\mathcal{N}_n}
	y(n',m)\frac{c_{n}^k(d)}{N^*_n} 
	\nonumber \\
	&+\frac{1}{q_n-B_n+B^{k*}_n(d)}\sum_{n'\in\mathcal{N}_n^{k*}(d)\cap\mathcal{B}_n}
	y(n',m)\frac{1}{N^*_n}
	\nonumber \\
	\ge & \frac{c_{n}^k(d)}{q_n-B_n+B^{k*}_n(d)} \!\!
	\sum_{n'\in\mathcal{N}_n^{k*}(d)}\!\!
	y(n',m)
	\frac{1}{N^*_n}\!\!
\end{align}
from which, plugging the definition \eqref{definition:ckn} of $c_{n}^k(d)$ and $|\mathcal{N}_n^{k*}(d)|=N_n+B_n-2q_n$ yields that
\begin{align}
	[\mathbf{Y}^k(d)]_{nm}
	\ge &\frac{q_n-B_n}{N_n+B_n-2q_n-B^{k*}_n(d)} \frac{N_n+B_n-2q_n}{q_n-B_n+B^{k*}_n(d)}
	\frac{1}{4N^*_n}
	\nonumber \\
	\ge &\frac{(q_n-B_n)(N_n+B_n-2q_n)}{(N_n-q_n)^2}
	\frac{1}{N^*_n},
\end{align}
where the last inequality comes from $4(N_n+B_n-2q_n-B^{k*}_n(d))(q_n-B_n+B^{k*}_n(d))\le(N_n-q_n)^2$.

As a result, in total $N_n-q_n+1$ elements have the lower bound $\mu_0$: $N_n+B_n-2q_n-B^{k*}_n(d)$ elements in $\mathcal{N}_n^{k*}(d)\cap\mathcal{N}_n$, $q_n-B_n+B^{k*}_n(d)$ elements in $\mathcal{L}_n^{k+}(d)$ or $\mathcal{L}_n^{k-}(d)$, as well as the element $n$.
Therefore Condition (c4) holds with constant
\begin{align}
	\mu_0 := \min_{n\in\mathcal{N}}
	&\left\{
	\frac{N_n+B_n-3q_n}{N_n-2q_n},
	\frac{(N_n+B_n-2q_n)(q_n-B_n)}{(N_n-q_n)^2}
	\right\}
	\times\frac{1}{N_n+B_n-2q_n+1}.
\end{align}

\noindent\textbf{Verification of (c5) for Case 2}. In Case 2, for $m\in \mathcal{N}_n^{k*}(d)\cap\mathcal{N}_n$, we have
\begin{align}
	[\mathbf{Y}^k(d)]_{nm} = \left(
	1-c_{n}^k(d)
	\right)\frac{1}{N^*_n}
	\le \frac{1}{N^*_n}.
\end{align}
For $m\in\mathcal{L}_n^{k+}(d)~\text{or}~\mathcal{L}_n^{k-}(d)$, $[\mathbf{Y}^k(d)]_{nm}$ reaches its maximum when all $y(n',m)=1$, such that
\begin{align}
	[\mathbf{Y}^k(d)]_{nm}
	=& \frac{1}{q_n-B_n+B^{k*}_n(d)}\sum_{n'\in\mathcal{N}_n^{k*}(d)\cap\mathcal{N}_n}
	y(n',m)\frac{1}{N^*_n} c_{n}^k(d) \nonumber \\
	&+\frac{1}{q_n-B_n+B^{k*}_n(d)}\sum_{n'\in\mathcal{N}_n^{k*}(d)\cap\mathcal{B}_n}
	y(n',m)\frac{1}{N^*_n}
	\nonumber \\
	\le & \frac{q_n-B_n}{q_n-B_n+B^{k*}_n(d)}\frac{1}{N^*_n}
	+\frac{B^{k*}_n(d)}{q_n-B_n+B^{k*}_n(d)}\frac{1}{N^*_n}
	= \frac{1}{N^*_n}. \nonumber
\end{align}

\noindent\textbf{Verification of (c6)}. Note that Assumption \ref{assumption:networkStructure} guarantees that any subgraph removing all agents in $\mathcal{W}$ and any additional $q_n$ incoming edges at honest agent $n$ has a source node $n^*$ and it can reach any other agent within $\tau_\mathcal{G}$ steps. Therefore, no element in the $n^*$-th column of matrix $(\mathbf{Y}^k(d))^{\tau_\mathcal{G}}$ is zero.

To show that there exist no more than $H_\mathcal{G}$ types of $\mathbf{Y}^k(d)$, notice that each $\mathbf{Y}^k(d)$ containing at least $N_n-q_n+1$ non-zero elements corresponds to a subgraph obtained by removing all Byzantine agents with their edges, and removing any additional $q_n$ incoming edges at honest node $n$. Every such subgraph is contained in $\mathcal{H}_\mathcal{G}$, whose cardinality is $H_\mathcal{G}$.

Now all conditions have been verified, and thus Theorems \ref{theorem:consensu-rate-theorem-global-time-varying} and \ref{theorem:convergence-trimmedMean-time-varying} can be applied to Byrd-TD($\lambda$).

\subsection{Proof of Corollary \ref{corollary:asymptoticTDerror}}

Recall the definition of $F(\bm{\theta})$ as
\begin{align}\label{ling-0002}
	&F(\bm{\theta}^{k+1}_n)
	=\frac{1}{2} \sum_{s\in\mathcal{S}} \rho(s)(\bm{\phi}(s)^{\top} \bm{\theta}^{k+1}_n - \mathcal{V}(s))^2
	\nonumber\\
	\le&\frac{1}{2} \sum_{s\in\mathcal{S}} \rho(s)
	[2(\bm{\phi}(s)^{\top} (\bm{\theta}^{k+1}_n - \bm{\theta}^{\infty}_{\lambda}))^2
	+2(\bm{\phi}(s)^{\top}\bm{\theta}^{\infty}_{\lambda} - \mathcal{V}(s))^2]
	\nonumber\\
	=&(\bm{\theta}^{k+1}_n - \bm{\theta}^{\infty}_{\lambda})^{\top} \mathbf{\Phi}_c (\bm{\theta}^{k+1}_n - \bm{\theta}^{\infty}_{\lambda})
	+2F(\bm{\theta}^{\infty}_{\lambda})
	\nonumber\\
	\le&\sigma_{c}\|\bm{\theta}^{k+1}_n - \bm{\theta}^{\infty}_{\lambda}\|^2
	+2F(\bm{\theta}^{\infty}_{\lambda}),
\end{align}
where the first inequality comes from the fact $(x+y)^2\le 2x^2+2y^2$, $\mathbf{\Phi}_{c}$ is defined by
\begin{align}
	\mathbf{\Phi}_{c}:= \mathbf{\Phi}^\top\mathbf{D}\mathbf{\Phi}\in\mathbb{R}^{D\times D},
\end{align}
and $\sigma_{c}$ is the largest singular value of $\mathbf{\Phi}_{c}$.

Combining with \eqref{inequality:approximationAccuracy}, \eqref{ling-0002} implies the following bound
\begin{align}
	\label{inequality:asymptoticTDerror-1}
	F(\bm{\theta}^{k+1}_n)
	\le& \sigma_{c}\|\bm{\theta}^{k+1}_n - \bm{\theta}^{\infty}_{\lambda}\|^2
	+2F(\bm{\theta}^{\infty}_{\lambda})
	\le \sigma_{c}\|\bm{\theta}^{k+1}_n - \bm{\theta}^{\infty}_{\lambda}\|^2
	+2\frac{1-\gamma\lambda}{1-\gamma}\min_{\bm{\theta}}F(\bm{\theta}).
\end{align}

Plugging \eqref{inequality:convergence_trimmedMean} into \eqref{inequality:asymptoticTDerror-1}, summing up over $n\in\mathcal{N}$ and taking limits on both sides of \eqref{inequality:asymptoticTDerror-1} yield
\begin{align}
	&\lim_{k\to\infty}\sup\frac{1}{N}\sum_{n\in\mathcal{N}}F(\bm{\theta}^{k+1}_n)
	\nonumber\\
	\le& \lim_{k\to\infty}\sup\frac{1}{N}\sum_{n\in\mathcal{N}}\sigma_{c}\|\bm{\theta}^{k+1}_n - \bm{\theta}^{\infty}_{\lambda}\|^2
	+2F(\bm{\theta}^{\infty}_{\lambda})
	\nonumber\\
	\le& 2C_4 \frac{D_{\mathcal{G}}\delta^2}{(1-\gamma\lambda)^2}
	+2\frac{1-\gamma\lambda}{1-\gamma}\min_{\bm{\theta}}F(\bm{\theta}).
\end{align}
This completes the proof.

\section{Proofs of supporting lemmas}
In this section, we will prove some lemmas used before.

\subsection{Proof of Lemma \ref{lemma:limit-of-product-of-Yk}}

Note that Condition (c6) guarantees that there exists at least one column of $(\mathbf{Y}^k(d))^{\tau_\mathcal{G}}$, denoted as the $n^*$-th column, being non-zero element-wise, and there exist no more than $H_\mathcal{G}$ choices of $\mathbf{Y}^k(d)$. Therefore, at least one $\mathbf{Y}^*(d)$ appears more than $\tau_\mathcal{G}$ times in the production $\bar{\mathbf{Y}}^{k_0+\tau_\mathcal{G}H_\mathcal{G}}(d, k_0)$ and there exists at least one column of $(\mathbf{Y}^k(d))^{\tau_\mathcal{G}}$, denoted as the $n^*$-th column, must be non-zero element-wise. In addition, since all non-zero elements in $\mathbf{Y}^k(d)$ are greater than $\mu_0$, all elements in the $n^*$-th column of $\bar{\mathbf{Y}}^{k_0+\tau_\mathcal{G}H_\mathcal{G}}(d, k_0)$ must be greater than $\mu_0^{\tau_\mathcal{G}H_\mathcal{G}}$.

Similar to \cite[Lemma 3]{vaidya2012matrix}, for a row-stochastic matrix $\mathbf{W}$, we define
\begin{align}
	\mathcal{U}_1(\mathbf{W}) :=
	1 - \min_{n_1, n_2\in\mathcal{N}} \sum_{m\in\mathcal{N}}\min\{
	[\mathbf{W}]_{n_1m},
	[\mathbf{W}]_{n_2m}
	\},
\end{align}
then we can derive the following bound using the technique similar to \cite[Lemma 4]{vaidya2012matrix}, as
\begin{align}
	\label{inequality:u1-of-Yk}
	\mathcal{U}_1(\bar{\mathbf{Y}}^{k_0+(S+1)\tau_\mathcal{G}H_\mathcal{G}-1}(d, k_0+S\tau_\mathcal{G}H_\mathcal{G}))
	\le 1-\mu_0^{\tau_\mathcal{G}H_\mathcal{G}},
\end{align}
In addition, we can define another function $\mathcal{U}_2(\mathbf{W})$ for a row-stochastic matrix $\mathbf{W}$ as
\begin{align}
	\mathcal{U}_2(\mathbf{W}) :=
	\max_{m\in\mathcal{N}}
	\max_{n_1, n_2\in\mathcal{N}}
	\left|[\mathbf{W}]_{n_1 m}-[\mathbf{W}]_{n_2 m}\right|.
\end{align}
It is easy to see that all rows of matrix $\mathbf{W}$ are identical if and only if $\mathcal{U}_1(\mathbf{W})=0$ or $\mathcal{U}_2(\mathbf{W})=0$.
Additionally, \cite{hajnal1958weak} has shown that $\mathcal{U}_1(\mathbf{W})$ and  $\mathcal{U}_2(\mathbf{W})$ satisfy
\begin{align}
	\mathcal{U}_2(\prod_{i=k_0}^{k}\mathbf{W}_i)
	\le \prod_{i=k_0}^{k}\mathcal{U}_1(\mathbf{W}_i).
\end{align}

We can see that $\mathcal{U}_1(\mathbf{W})\le1$ and further conclude that
\begin{align}
	\mathcal{U}_2(\bar{\mathbf{Y}}_{k_0}^k(d))
	\le &
	\mathcal{U}_1\left(
	\bar{\mathbf{Y}}^{k}\left(d, k_0+\lfloor\frac{k-k_0+1}{\tau_\mathcal{G}H_\mathcal{G}}\rfloor\tau_\mathcal{G}H_\mathcal{G}\right)
	\right)
	\prod_{s=0}^{\lfloor\frac{k-k_0+1}{\tau_\mathcal{G}H_\mathcal{G}}\rfloor-1} \mathcal{U}_1(\bar{\mathbf{Y}}^{k_0+(s+1)\tau_\mathcal{G}H_\mathcal{G}-1}(d, k_0+s\tau_\mathcal{G}H_\mathcal{G}))
	\nonumber \\
	\le &
	\prod_{s=0}^{\lfloor\frac{k-k_0+1}{\tau_\mathcal{G}H_\mathcal{G}}\rfloor-1} \mathcal{U}_1(\bar{\mathbf{Y}}^{k_0+(s+1)\tau_\mathcal{G}H_\mathcal{G}-1}(d, k_0+s\tau_\mathcal{G}H_\mathcal{G}))
	\le (1-\mu_0^{\tau_\mathcal{G}H_\mathcal{G}})^{\lfloor\frac{k-k_0+1}{\tau_\mathcal{G}H_\mathcal{G}}\rfloor}
	\nonumber \\
	\le& (1-\mu_0^{\tau_\mathcal{G}H_\mathcal{G}})^{\frac{k-k_0+1}{\tau_\mathcal{G}H_\mathcal{G}}-1},
\end{align}
where the third inequality comes from inequality \eqref{inequality:u1-of-Yk}.

From Condition (c3), at least $N_n-q_n+1$ elements in each row are lower bounded by $\mu_0$, implying that $\mu_0\le\frac{1}{N_n-q_n+1}\le\frac{1}{2}$ and  $(1-\mu_0^{\tau_\mathcal{G}H_\mathcal{G}})^{-1}\le 2$. This relation can further bound $\mathcal{U}_2(\bar{\mathbf{Y}}_{k_0}^k(d))$ by
\begin{align}
	\mathcal{U}_2(\bar{\mathbf{Y}}_{k_0}^k(d))
	\le (1-\mu_0^{\tau_\mathcal{G}H_\mathcal{G}})^{\frac{k-k_0+1}{\tau_\mathcal{G}H_\mathcal{G}}-1}
	\le 2(1-\mu_0^{\tau_\mathcal{G}H_\mathcal{G}})^{\frac{k-k_0+1}{\tau_\mathcal{G}H_\mathcal{G}}}.
\end{align}
With this inequality, we can show that for any possible $\{\mathbf{Y}^k(d)\}$, the limit of $\bar{\mathbf{Y}}_{k_0}^k(d)$ exists. Recalling that $\bar{\mathbf{Y}}_{k_0}^{k+1}(d)=\mathbf{Y}^{k+1}(d)\bar{\mathbf{Y}}_{k_0}^k(d)$ and $\mathbf{Y}^{k+1}(d)$ is a row-stochastic matrix, it holds for all $m\in\mathcal{N}$ that
\begin{align}
	\min_{n'}[\bar{\mathbf{Y}}_{k_0}^{k+1}(d)]_{n'm}
	\ge& \min_{n'}[\bar{\mathbf{Y}}_{k_0}^k(d)]_{n'm},
	\nonumber \\
	\max_{n'}[\bar{\mathbf{Y}}_{k_0}^{k+1}(d)]_{n'm}
	\le& \max_{n'}[\bar{\mathbf{Y}}_{k_0}^k(d)]_{n'm}
\end{align}
Together with inequality
\begin{align}
	\max_{n'}[\bar{\mathbf{Y}}_{k_0}^k(d)]_{n'm}
	-\min_{n'}[\bar{\mathbf{Y}}_{k_0}^k(d)]_{n'm}
	\le \mathcal{U}_2(\bar{\mathbf{Y}}_{k_0}^k(d))
	\to 0
	~~\text{as}~~k\to+\infty,
\end{align}
it reaches that the limit of $\bar{\mathbf{Y}}_{k_0}^k(d)$ exists and all elements of this limit in the same column share the same value. Hence, we use $\bm{1}\bm{p}_{k_0}^{\top}(d)$ to denote this limit, namely,
\begin{align}
	\lim_{k\to\infty}\bar{\mathbf{Y}}_{k_0}^k(d) &= \bm{1}\bm{p}_{k_0}^{\top}(d).
	\label{inequality:consensus-rate-norm-1}
\end{align}
In addition, we can also see that $\bar{\mathbf{Y}}_{k_0}^k(d)$ converges to its limit at a linear rate. To see this fact, we note that for any $m\in\mathcal{N}$
\begin{align}
	\label{inequality:p-bound-bar}
	\min_{n\in\mathcal{N}} [\bar{\mathbf{Y}}_{k_0}^k(d)]_{nm}
	\le [\bm{p}_{k_0}(d)]_m \le \max_{n\in\mathcal{N}} [\bar{\mathbf{Y}}_{k_0}^k(d)]_{nm},
\end{align}
and hence for any $n, m\in\mathcal{N}$, we have
\begin{align}
	\left|[\bar{\mathbf{Y}}_{k_0}^k(d)]_{nm}-[\bm{p}_{k_0}(d)]_m\right|
	\le& \max_{n_1, n_2\in\mathcal{N}}
	\left|
	[\bar{\mathbf{Y}}_{k_0}^k(d)]_{n_1m}
	-[\bar{\mathbf{Y}}_{k_0}^k(d)]_{n_2m}
	\right|
	\nonumber \\
	\le& \max_{m'\in\mathcal{N}}
	\max_{n_1, n_2\in\mathcal{N}}
	\left|
	[\bar{\mathbf{Y}}_{k_0}^k(d)]_{n_1m'}
	-[\bar{\mathbf{Y}}_{k_0}^k(d)]_{n_2m'}
	\right|
	\nonumber \\
	=& \mathcal{U}_2(\bar{\mathbf{Y}}_{k_0}^k(d))
	\le 2(1-\mu_0^{\tau_\mathcal{G}H_\mathcal{G}})^{\frac{k-k_0+1}{\tau_\mathcal{G}H_\mathcal{G}}}.
\end{align}
Therefore, the distance $\|\bar{\mathbf{Y}}_{k_0}^k(d) - \bm{1}\bm{p}_{k_0}^{\top}(d)\|^2$ has the bound that
\begin{align}
	\|\bar{\mathbf{Y}}_{k_0}^k(d) - \bm{1}\bm{p}_{k_0}^{\top}(d)\|^2 
	\le& \|\bar{\mathbf{Y}}_{k_0}^k(d) - \bm{1}\bm{p}_{k_0}^{\top}(d)\|^2_F
	\nonumber \\
	=& \sum_{n,m\in\mathcal{N}} ([\bar{\mathbf{Y}}_{k_0}^k(d)]_{nm}-[\bm{p}_{k_0}(d)]_m)^2
	\nonumber \\
	\leq& 4N^2(1-\mu_0^{\tau_\mathcal{G}H_\mathcal{G}})^{\frac{2(k-k_0+1)}{\tau_\mathcal{G}H_\mathcal{G}}},
	\label{inequality:consensus-rate-norm-2}
\end{align}
which is inequality \eqref{equation:lemma-Y-bar-limit-rate} with notation $\mu$ defined as
\begin{equation}
	\label{definition:mu}
	\mu := (1-\mu_0^{\tau_\mathcal{G}H_\mathcal{G}})^{\frac{2}{\tau_\mathcal{G}H_\mathcal{G}}}.
\end{equation}


Similar to inequality \eqref{inequality:p-bound-bar}, since all $\{\mathbf{Y}^k(d)\}$ are row-stochastic matrices, the elements of $\bm{p}_{k_0}(d)$ are bounded by the corresponding columns of $\mathbf{Y}^{k_0}(d)$, such that
\begin{align}
	\label{inequality:p-bound}
	\min_{n\in\mathcal{N}} [\mathbf{Y}^{k_0}(d)]_{nm}
	\le [\bm{p}_{k_0}(d)]_m \le \max_{n\in\mathcal{N}} [\mathbf{Y}^{k_0}(d)]_{nm}.
\end{align}
From Conditions (c1) and (c5), it follows that 
\begin{equation}
	0\le[\mathbf{Y}^{k_0}(d)]_{nm}\le \frac{1}{\min_{n\in\mathcal{N}}\{N_n+B_n-2q_n+1\}}.
\end{equation}
Therefore, the elements of $\bm{p}_{k_0}(d)$ are also bounded in $\left[0, \frac{1}{\min_{n\in\mathcal{N}}\{N_n+B_n-2q_n+1\}}\right]$. Now we have completed the proof.

\subsection{Proof of Lemma \ref{lemma:different-local-global-update}}

From the definition in \eqref{definition:eligibilityTrace}, $\bm{z}^k$ has a bounded norm
\begin{align}
	\|\bm{z}^k\|
	= \Big\|\sum_{\kappa=0}^k (\gamma\lambda)^{k-\kappa}\bm{\phi}(s^\kappa)\Big\|
	\le \sum_{\kappa=0}^k (\gamma\lambda)^{k-\kappa}\|\bm{\phi}(s^\kappa)\|
	\le \sum_{\kappa=0}^k (\gamma\lambda)^{k-\kappa}
	\le \frac{1}{1-\gamma\lambda},
\end{align}
where the second inequality is the result of normalized features in Assumption \ref{assumption:NormalizedFeature}. Observe from the definitions in \eqref{definition:Akbk} that $\mathbf{A}^k$ has a bounded norm as well, given by
\begin{align}
	\label{inequality:Ak-bound}
	\|\mathbf{A}^k\|
	\le& \|\bm{z}^k(\gamma\bm{\phi}(s^{k+1})-\bm{\phi}(s^k))\|
	\le \|\bm{z}^k\|\|\gamma\bm{\phi}(s^{k+1})-\bm{\phi}(s^k)\|
	\nonumber\\
	\le& \|\bm{z}^k\|(\gamma\|\bm{\phi}(s^{k+1})\|+\|\bm{\phi}(s^k)\|)
	= \frac{1+\gamma}{1-\gamma\lambda},
\end{align}
and $b^k_n$ has a bounded variation, given by
\begin{align}
	\frac{1}{N}\sum_{n\in\mathcal{N}}\|\bm{b}^k_n-\bar{\bm{b}}^k\|^2
	= &\frac{1}{N}\sum_{n\in\mathcal{N}}|r_n^k-\bar r^k|^2\|\bm{z}^k\|^2
	\nonumber\\
	\le& \frac{1}{N}\sum_{n\in\mathcal{N}}|r_n^k-\bar r_n^k|^2\frac{1}{(1-\gamma\lambda)^2}
	\le \frac{\delta^2}{(1-\gamma\lambda)^2},
\end{align}
where the last inequality follows from the bounded reward variation in Assumption \ref{assumption:rewardVariation}.

With these inequalities, we bound the difference between the local and global increments by
\begin{align}
	\frac{1}{N}\sum_{n\in\mathcal{N}}\|\bm{g}^k_n(\bm{\theta}_n)-\bm{\bar {g}}^k(\bm{\theta}')\|^2
	=&\frac{1}{N}\sum_{n\in\mathcal{N}}\|\mathbf{A}^k(\bm{\theta}_n-\bm{\theta}')+
	\bm{b}^k_n-\bar{\bm{b}}^k\|^2
	\nonumber \\
	\le&\frac{1}{N}\sum_{n\in\mathcal{N}}2\|\mathbf{A}^k\|^2\|\bm{\theta}_n-\bm{\theta}'\|^2
	+\frac{1}{N}\sum_{n\in\mathcal{N}}2\|\bm{b}^k_n-\bar{\bm{b}}^k\|^2
	\nonumber \\
	\le&2\left(\frac{1+\gamma}{1-\gamma\lambda}\right)^2\frac{1}{N}\sum_{n\in\mathcal{N}}\|\bm{\theta}_n-\bm{\theta}'\|^2
	+ \frac{2\delta^2}{(1-\gamma\lambda)^2},
\end{align}
which completes the proof.

\subsection{Proof of Lemma \ref{lemma:convergence-omega}}
We use $k_m$ to represent the smallest $k$ satisfying $k>\tau(\eta^k)$. When $k>k_m$, taking expectation on both sides and substituting \eqref{inequality:convergence-rate-Central-TD-itretative} in Theorem \ref{theorem:convergence-rate-Central-TD-itretative} into \eqref{inequality:convergence-rate-consensus-mapping-1.5}, it follows
\begin{align}
	\label{inequality:convergence-rate-consensus-mapping-2}
	\mathbb{E}\|\bm{\omega}^{k+1}-\bm{\theta^{\infty}}_{\lambda}\|^2
	\le \frac{1-\sigma_{\rm min}\eta^k}{1-\beta^k}
	\mathbb{E}\|\bm{\omega}^{k}-\bm{\theta^{\infty}}_{\lambda}\|^2
	+\frac{1}{1-\beta^k} 2C_6 (\eta^k)^2 C \ln(\frac{k+k^0}{\eta})
	+\frac{(\eta^k)^2}{\beta^k} \frac{2D_{\mathcal{G}}\delta^2}{(1-\gamma\lambda)^2}.
\end{align}
With $\beta^k=\sigma_{\rm min}\eta^k/2$ and $\sigma_{\rm min}\eta^0/2\le 1/2$, it follows that
$\frac{1-\sigma_{\rm min}\eta^k}{1-\beta^k}    \le 1-\frac{\sigma_{\rm min}\eta^k}{2}$ and $\frac{1}{1-\beta^k} \le 2$,
with which \eqref{inequality:convergence-rate-consensus-mapping-2} becomes
\begin{align}
	\label{inequality:convergence-rate-consensus-mapping-3}
	\mathbb{E}\|\bm{\omega}^{k+1}-\bm{\theta^{\infty}}_{\lambda}\|^2
	\le& (1-\frac{\sigma_{\rm min}\eta^k}{2})\mathbb{E}\|\bm{\omega}^{k}-\bm{\theta^{\infty}}_{\lambda}\|^2
	+4C_6 (\eta^k)^2 C \ln(\frac{k+k^0}{\eta})
	+\frac{2\eta^k}{\sigma_{\rm min}}  \frac{2D_{\mathcal{G}}\delta^2}{(1-\gamma\lambda)^2} \nonumber\\
	\le& \prod_{j=k_m}^k(1-\frac{\sigma_{\rm min}\eta^j}{2})\mathbb{E}\|\bm{\omega}^{k_m}-\bm{\theta^{\infty}}_{\lambda}\|^2
	+\sum_{i=k_m}^{k}4C_6 (\eta^i)^2 C \ln(\frac{i+k^0}{\eta}) \prod_{j=i+1}^k(1-\frac{\sigma_{\rm min}\eta^j}{2})\nonumber \\
	&+\sum_{i=k_m}^{k}\frac{2\eta^i}{\sigma_{\rm min}}  \frac{2D_{\mathcal{G}}\delta^2}{(1-\gamma\lambda)^2} \prod_{j=i+1}^k(1-\frac{\sigma_{\rm min}\eta^j}{2}) \nonumber\\
	\le& \prod_{j=k_m}^k(1-\frac{\sigma_{\rm min}\eta^j}{2})\mathbb{E}\|\bm{\omega}^{k_m}-\bm{\theta^{\infty}}_{\lambda}\|^2
	+\sum_{i=0}^{k}4C_6 (\eta^i)^2 C \ln(\frac{i+k^0}{\eta}) \prod_{j=i+1}^k(1-\frac{\sigma_{\rm min}\eta^j}{2})\nonumber \\
	&+\sum_{i=0}^{k}\frac{2\eta^i}{\sigma_{\rm min}}  \frac{2D_{\mathcal{G}}\delta^2}{(1-\gamma\lambda)^2} \prod_{j=i+1}^k(1-\frac{\sigma_{\rm min}\eta^j}{2}).
\end{align}

To obtain a more compact form, consider the relation
\begin{align}
	1-\frac{\sigma_{\rm min}\eta^k}{2}
	= 1-\frac{\sigma_{\rm min}\eta}{2(k+k^0)}
	\le (1-\frac{1}{k+k^0})^{\sigma_{\rm min}\eta/2}
	= (\frac{k-1+k^0}{k+k^0})^{\epsilon}.
\end{align}
Thus, it follows
\begin{align}
	&\prod_{j=i+1}^k(1-\frac{\sigma_{\rm min}\eta^j}{2})
	\le(\frac{i+k^0}{k+k^0})^{\epsilon}.
\end{align}

This leads to the following bound on the second term at the RHS of \eqref{inequality:convergence-rate-consensus-mapping-3}, given by
\begin{align}
	\label{inequality:convergence-rate-consensus-mapping-3-1}
	\sum_{i=0}^{k}4C_6 (\eta^i)^2 C \ln(\frac{i+k^0}{\eta}) \prod_{j=i+1}^k(1-\frac{\sigma_{\rm min}\eta^j}{2})
	\le  &\sum_{i=0}^{k}4C_6 (\eta^i)^2 C \ln(\frac{i+k^0}{\eta}) (\frac{i+k^0}{k+k^0})^{\epsilon} \nonumber\\
	= & 4C_6C(\frac{\eta}{k+k^0})^{\epsilon}  \sum_{i=0}^{k} \frac{\ln(\frac{i+k^0}{\eta})}{(\frac{i+k^0}{\eta})^{2-\epsilon}}
	\le  4C_6C \varphi_\epsilon(k),
\end{align}
where the last inequality comes from the following lemma.
\begin{Lemma}
	\label{lemma:dominantFunction}
	When $k^0$ is chosen large enough that $k^0 > \eta e^{2-\epsilon}$, then
	\begin{align}
		(\frac{\eta}{k+k^0})^{\epsilon}\sum_{i=0}^{k} \frac{\ln(\frac{i+k^0}{\eta})}{(\frac{i+k^0}{\eta})^{2-\epsilon}}
		\le \varphi_\epsilon(k),
	\end{align}
	where the function $\varphi_\epsilon(k)$ is defined in \eqref{definition:convergence-rate-dominant}.
\end{Lemma}

Similarly, the third term at the RHS of \eqref{inequality:convergence-rate-consensus-mapping-3} satisfies
\begin{align}
	\label{inequality:convergence-rate-consensus-mapping-3-term3}
	\sum_{i=0}^{k}\frac{2\eta^i}{\sigma_{\rm min}} \frac{2D_{\mathcal{G}}\delta^2}{(1-\gamma\lambda)^2} \prod_{j=i+1}^k(1-\frac{\sigma_{\rm min}\eta^j}{2}) \le  &\sum_{i=0}^{k}\frac{2\eta}{(i+k^0)\sigma_{\rm min}} \frac{2D_{\mathcal{G}}\delta^2}{(1-\gamma\lambda)^2} (\frac{i+k^0}{k+k^0})^{\epsilon} \nonumber\\
	\le& \frac{2\eta}{\sigma_{\rm min}} \frac{2D_{\mathcal{G}}\delta^2}{(1-\gamma\lambda)^2}
	\sum_{i=0}^{k}(\frac{i+k^0}{k+k^0})^{\epsilon-1}  \frac{1}{k+k^0}
	\nonumber\\
	=&     \frac{2}{\sigma_{\rm min}^2} \frac{2D_{\mathcal{G}}\delta^2}{(1-\gamma\lambda)^2}
	\epsilon\sum_{i=0}^{k}(\frac{i+k^0}{k+k^0})^{\epsilon-1} \frac{1}{k+k^0}\le \frac{8}{\sigma_{\rm min}} \frac{D_{\mathcal{G}}\delta^2}{(1-\gamma\lambda)^2},
\end{align}
where the last inequality is the result of the following lemma.
\begin{Lemma}
	\label{lemma:asymptotic-error-summation}
	For any positive integer $k$ and large enough $k^0$, it holds
	\begin{align}
		\epsilon\sum_{i=0}^{k}(\frac{i+k^0}{k+k^0})^{\epsilon-1}\le 2.
	\end{align}
\end{Lemma}

Substituting \eqref{inequality:convergence-rate-consensus-mapping-3-term3} into \eqref{inequality:convergence-rate-consensus-mapping-3} yields
\begin{align}
	\label{inequality:convergence-rate-consensus-mapping-4}
	\mathbb{E}\|\bm{\omega}^{k+1}-\bm{\theta^{\infty}}_{\lambda}\|^2
	\le& \mathbb{E}\|\bm{\omega}^{k_m}-\bm{\theta^{\infty}}_{\lambda}\|^2\left(\frac{k_m+k^0}{k+k^0}\right)^\epsilon
	+4C_6C\varphi_{\epsilon}(k)
	+\frac{8}{\sigma_{\rm min}^2} \frac{D_{\mathcal{G}}\delta^2}{(1-\gamma\lambda)^2}.
\end{align}

Choosing $\beta^k=\frac{1}{2}$ in  \eqref{inequality:convergence-rate-consensus-mapping-1.5} and using Theorem \ref{theorem:convergence-rate-Central-TD-itretative}, we have
\begin{align}
	\label{inequality:convergence-rate-consensus-mapping-km-to-0}
	\|\bm{\omega}^{k_m}-\bm{\theta^{\infty}}_{\lambda}\|^2
	\le &2\|\bm{\omega}^{k_m-1}
	+\eta^k\bm{\bar g}^k(\bm{\omega}^{k})
	-\bm{\theta^{\infty}}_{\lambda}\|^2
	+2(\eta^k)^2\frac{2D_{\mathcal{G}}\delta^2}{(1-\gamma\lambda)^2}
	\nonumber \\
	\le& 2\|\bm{\omega}^{k_m-1}
	+\eta^k\bm{\bar g}^k(\bm{\omega}^{k})
	-\bm{\theta^{\infty}}_{\lambda}\|^2
	+2(\eta^k)^2\frac{2D_{\mathcal{G}}\delta^2}{(1-\gamma\lambda)^2}
	\nonumber \\
	\le& 4\left(1+\eta^k\left(\frac{1+\gamma}{1-\gamma\lambda}\right)\right)^2
	\|\bm{\omega}^{k_m-1}-\bm{\theta}^{\infty}_{\lambda}\|^2
	+2(\eta^k)^2C_7
	\nonumber \\
	\le& \left(2+2\eta^0\left(\frac{1+\gamma}{1-\gamma\lambda}\right)\right)^{2k_m}
	\|\bm{\omega}^{0}-\bm{\theta}^{\infty}_{\lambda}\|^2
	+2(\eta^0)^2C_7C_8,
\end{align}
where constants $C_7$ and $C_8$ are given by
\begin{align}
	C_7 :=& 2\left(\left(\frac{1+\gamma}{1-\gamma\lambda}\right)\| \bm{\theta}^{\infty}_{\lambda}\|+\left(\frac{R_{max}}{1-\gamma\lambda}\right)\right)^2
	+\frac{2D_{\mathcal{G}}\delta^2}{(1-\gamma\lambda)^2},
	\nonumber \\
	C_8 :=& \frac{\left(2+2\eta^0\left(\frac{1+\gamma}{1-\gamma\lambda}\right)\right)^{2k_m}-1}{\left(2+2\eta^0\left(\frac{1+\gamma}{1-\gamma\lambda}\right)\right)-1}.
\end{align}

Plugging \eqref{inequality:convergence-rate-consensus-mapping-km-to-0} into \eqref{inequality:convergence-rate-consensus-mapping-4}, we have
\begin{align}
	\label{inequality:convergence-rate-consensus-mapping-5}
	\mathbb{E}\|\bm{\omega}^{k}-\bm{\theta^{\infty}}_{\lambda}\|^2
	\le\frac{C_2}{2}
	\left(\frac{1}{k+k^0}\right)^\epsilon
	+\frac{C_3}{2}\varphi_{\epsilon}(k)
	+\frac{C_4}{2} \frac{D_{\mathcal{G}}\delta^2}{(1-\gamma\lambda)^2},
\end{align}
where we use a the fact that $\bm{\omega}^0=\bm{\theta}^{0}$ to simplify the inequality and the constant is defined as
\begin{align}
	C_2 := 2\|\bm{\theta^0}-\bm{\theta^{\infty}}_{\lambda}\|^2&
	\left(
	\left(2+2\eta^0\left(\frac{1+\gamma}{1-\gamma\lambda}\right)\right)^{2k_m}
	\|\bm{\theta}^{0}-\bm{\theta}^{\infty}_{\lambda}\|^2
	+2(\eta^0)^2C_7C_8
	\right)
	(k_m+k^0)^\epsilon,
	\nonumber \\
	&C_3 := 8C_6C
	\quad\quad\text{and}\quad\quad
	C_4 := \frac{16}{\sigma_{\rm min}^2}.
\end{align}

\subsection{Proof of Lemma \ref{lemma:dominantFunction}}

We proceed by discussing the bound under different values of $\epsilon$.

First, when $\epsilon\ge2$, $ \frac{\ln(x)}{x^{2-\epsilon}}$ is a monotonically increasing function, then
\begin{align*}
	\sum_{i=0}^{k} \frac{\ln(\frac{i+k^0}{\eta})}{(\frac{i+k^0}{\eta})^{2-\epsilon}}
	\le &\int_{\frac{k^0}{\eta}}^{\frac{k+1+k^0}{\eta}} \frac{\ln(x)}{x^{2-\epsilon}} dx  
	\le  \frac{1}{(\epsilon-1)^2}
	\left((\epsilon-1)\ln\left(\frac{k+1+k^0}{\eta}\right)+1\right)
	\left(\frac{k+1+k^0}{\eta}\right)^{\epsilon-1}.
\end{align*}

Second, when $1<\epsilon<2$ and $\frac{k^0}{\eta}\ge e^{2-\epsilon}$, $ \frac{\ln(x)}{x^{2-\epsilon}}$ is monotonically decreasing on $(\frac{k^0-1}{\eta}, \frac{k+k^0}{\eta})$, then
\begin{align*}
	\sum_{i=0}^{k} \frac{\ln(\frac{i+k^0}{\eta})}{(\frac{i+k^0}{\eta})^{2-\epsilon}}
	\le \int_{\frac{k^0-1}{\eta}}^{\frac{k+k^0}{\eta}} \frac{\ln(x)}{x^{2-\epsilon}} dx
	\le \frac{1}{(\epsilon-1)^2}
	\left((\epsilon-1)\ln\left(\frac{k+k^0}{\eta}\right)+1\right)
	\left(\frac{k+k^0}{\eta}\right)^{\epsilon-1}.
\end{align*}

Third, when $\epsilon=1$ and $\frac{k^0}{\eta}\ge e^{2-\epsilon}$, $ \frac{\ln(x)}{x^{2-\epsilon}}$ is monotonically decreasing on $(\frac{k^0-1}{\eta}, \frac{k+k^0}{\eta})$, then
\begin{align*}
	\sum_{i=0}^{k} \frac{\ln(\frac{i+k^0}{\eta})}{(\frac{i+k^0}{\eta})^{2-\epsilon}}
	\le \int_{\frac{k^0-1}{\eta}}^{\frac{k+k^0}{\eta}} \frac{\ln(x)}{x} dx
	\le \ln\left(\frac{k+k^0}{\eta}\right)^2.
\end{align*}

Fourth, when $0<\epsilon<1$ and $\frac{k^0}{\eta}\ge e^{2-\epsilon}$, $ \frac{\ln(x)}{x^{2-\epsilon}}$ is monotonically decreasing on $(\frac{k^0-1}{\eta}, \frac{k+k^0}{\eta})$, then
\begin{align*}
	\sum_{i=0}^{k} \frac{\ln(\frac{i+k^0}{\eta})}{(\frac{i+k^0}{\eta})^{2-\epsilon}}
	\le& \int_{\frac{k^0-1}{\eta}}^{\frac{k+k^0}{\eta}} \frac{\ln(x)}{x^{2-\epsilon}} dx
	\le \frac{1}{(\epsilon-1)^2}
	\left((1-\epsilon)\ln\left(\frac{k^0-1}{\eta}\right)-1\right)
	\left(\frac{k^0-1}{\eta}\right)^{\epsilon-1}.
\end{align*}

For simplicity, we merge the first two cases and complete the proof.

\subsection{Proof of Lemma \ref{lemma:asymptotic-error-summation}}

When $\epsilon>1$, $x^{\epsilon-1}$ is monotonically increasing on $(0,+\infty)$ and
\begin{align}
	\sum_{i=0}^{k}(\frac{i+k^0}{k+k^0})^{\epsilon-1} \frac{1}{k+k^0}
	=&\sum_{i=0}^{k-1}(\frac{i+k^0}{k+k^0})^{\epsilon-1} \frac{1}{k+k^0}+\frac{1}{k+k^0}
	\nonumber\\
	\le& \int_{0}^{1}x^{\epsilon-1}dx +\frac{1}{k+k^0}
	\le \frac{1}{\epsilon}+\frac{1}{k+k^0}.
\end{align}

When $\epsilon=1$, it holds
\begin{align}
	\sum_{i=0}^{k}(\frac{i+k^0}{k+k^0})^{\epsilon-1} \frac{1}{k+k^0}
	= \frac{k+1}{k+k^0}
	\le 1.
\end{align}

When $0<\epsilon<1$, $x^{\epsilon-1}$ is monotonically decreasing on $(0,+\infty)$ and
\begin{align}
	\sum_{i=0}^{k}(\frac{i+k^0}{k+k^0})^{\epsilon-1} \frac{1}{k+k^0}
	\le \int_{0}^{1}x^{\epsilon-1}dx
	\le \frac{1}{\epsilon}.
\end{align}

In each of the cases above, $\epsilon\sum_{i=0}^{k}(\frac{i+k^0}{k+k^0})^{\epsilon-1}\le 2$, which completes the proof.

\end{document}